\documentclass[onefignum,onetabnum]{siamart250211}

\usepackage{amsfonts, amssymb, amsopn}
\usepackage[mathscr]{eucal}
\usepackage{bm, upgreek}
\usepackage{graphicx, subcaption, epstopdf, booktabs}
\usepackage[linesnumbered,ruled,vlined,algo2e]{algorithm2e}
\SetAlFnt{\small}

\ifpdf
  \DeclareGraphicsExtensions{.eps,.pdf,.png,.jpg}
\else
  \DeclareGraphicsExtensions{.eps}
\fi

\renewenvironment{equation*}{\[}{\]\ignorespacesafterend}

\newsiamremark{remark}{Remark}
\crefname{algocf}{algorithm}{algorithms}
\Crefname{algocf}{Algorithm}{Algorithms}

\headers{Matrices over a Hilbert space}{S. Budzinskiy}

\title{Matrices over a Hilbert space and their low-rank cross approximation\thanks{This research was funded by the Disruptive Innovation – Early Career Seed Money funding program of the Austrian Academy of Sciences (\"OAW) and the Austrian Science Fund (FWF).}}

\author{
Stanislav Budzinskiy\thanks{Faculty of Mathematics, University of Vienna, Kolingasse 14-16, 1090 Vienna, Austria (\email{stanislav.budzinskiy\allowbreak@univie.ac.at}).}
}

\newcommand{\Clx}{\mathbb{C}}
\newcommand{\Real}{\mathbb{R}}
\newcommand{\F}{\mathbb{F}}
\newcommand{\N}{\mathbb{N}}

\newcommand{\Hilb}{\mathrm{H}}
\newcommand{\Ell}{\mathrm{L}}

\newcommand{\set}[3][]{#1\{ #2 : #3 #1\}}
\newcommand{\argmax}{\mathrm{argmax}}

\newcommand{\Bensor}[1]{\bm{\Tensor{#1}}}
\newcommand{\Tensor}[1]{{\mathscr{#1}}}
\newcommand{\Batrix}[1]{\bm{\Matrix{#1}}}
\newcommand{\Matrix}[1]{\mathsf{#1}}
\newcommand{\Bentry}[1]{{\bm{#1}}}
\newcommand{\trans}{\intercal}
\newcommand{\pinv}{\dagger}
\newcommand{\trace}{\mathrm{tr}}

\newcommand{\krp}{\mathbin{\otimes}}
\newcommand{\rank}[2][]{\mathrm{rank} #1( #2 #1)}
\newcommand{\Rank}[3][]{\mathrm{rank_{#2}} #1( #3 #1)}
\newcommand{\truncate}[3][]{#1\lfloor #2 #1\rfloor_{#3}}
\newcommand{\Index}[1]{{#1}}
\newcommand{\mIndex}[1]{\mathfrak{#1}}
\newcommand{\cross}[3][]{\mathrm{cross}#1(#2, #3 #1)}

\newcommand{\Norm}[3][]{#1\| #2 #1\|_{\mathrm{#3}}}

\newcommand{\Dotp}[4][]{#1\langle #2, #3 #1\rangle_{\mathrm{#4}}}

\newcommand{\scup}[1]{\textup{\textsc{#1}}}

\begin{document}
\maketitle

\begin{abstract}
Motivated by applications in reduced-order modeling (ROM) of parametric partial differential equations, we investigate the algebraic properties of Bochner matrices---matrices with entries in an abstract Hilbert space.
Low-rank cross approximation is extended to Bochner matrices and its approximation guarantees are derived.
The high-dimensional nature of the entries is shown to manifest itself in maximum-volume bounds, making them smaller than in the classical setting.
An analogue of adaptive cross approximation is proposed and validated as a non-intrusive ROM method in numerical experiments with parametric nonlinear Stokes equations.
\end{abstract}

\begin{keywords}
vector-valued matrices, low-rank approximation, cross approximation, maximum volume, reduced-order modeling, parametric partial differential equations
\end{keywords}

\begin{MSCcodes}
15A69, 65D40, 65F55, 65N99 
\end{MSCcodes}

\section{Introduction}
The kind of entries a matrix has is determined by an application from which it arises.
The most common are real and complex matrices, appearing universally across scientific computing, signal processing, and data science.
Matrices over finite fields \cite[Ch.~7]{hachenberger2020matrices} are used in coding theory \cite{silva2010communication} and cryptography \cite{kleinjung2010factorization}.

Matrices over commutative rings \cite{brown1992matrices} have also been extensively studied.
Matrices of polynomials \cite{foster2009algorithm}, Laurent polynomials \cite{cescato2010qr}, and analytic functions \cite{bunse1991numerical, weiss2024properties} appear in multichannel broadband signal processing \cite{neo2023polynomial}.
Matrices over the convolution ring of vectors \cite{kilmer2011factorization, kilmer2013third, mor2025quasitubal} are well-suited to represent visual data \cite{dian2019hyperspectral}.

Matrices over the non-commutative ring of quaternions \cite{zhang1997quaternions} are used in signal and image processing \cite{miron2006quaternion, jia2019robust, miron2023quaternions}, while matrices of linear operators appear in the context of partial differential equations (PDEs) \cite{tretter2008spectral, barbarino2020block}.

The motivating application behind the present work, however, produces matrices whose entries form neither a field nor a ring, but rather reside in a Hilbert space.

\subsection{Reduced-order modeling of parametric PDEs}
\label{sec:intro:rom}
Let $\Omega \subset \Real^{N}$ be open and bounded, and consider a Dirichlet boundary-value problem for the Poisson equation in $\Omega$ with right-hand side $f \in \Ell^2(\Omega)$ and varying coefficient $\kappa \in \Ell^\infty(\Omega)$ such that $\kappa \geq 1$ almost everywhere:
\begin{equation*}
    -\mathrm{div}(\kappa \nabla u) = f~\text{ in }~\Omega, \qquad u = 0~\text{ on }~\partial\Omega.
\end{equation*}
This problem admits a unique weak solution in the Sobolev space $\Hilb^{1}_0(\Omega)$ \cite[Ch.~9]{brezis2011functional}.
Suppose now that $f = f_{\alpha,\beta}$ and $\kappa = \kappa_{\alpha,\beta}$ depend on parameters $\alpha,\beta \in [0,1]$; then there is a well-defined \emph{solution map} $(\alpha, \beta) \mapsto u_{\alpha,\beta}$.
Upon introducing grids $\{ \alpha_i \}_{i = 1}^{m}$ and $\{ \beta_j \}_{j = 1}^{n}$, we obtain an $m \times n$ matrix with entries in $\Hilb^{1}_0(\Omega)$.

In reduced-order modeling (ROM), one often seeks to construct a data-driven surrogate model to approximate a quantity of interest (QoI) associated with the solution of a parametric PDE, allowing for rapid evaluation without querying the PDE solver.
This essentially reduces to approximating a function $[0,1]^2 \to \Real$ from samples, a problem with numerous established approaches.
When the QoI is discretized into a matrix, \emph{cross approximation} \cite{goreinov1997theory, bebendorf2000approximation} provides a surrogate model by inspecting a small number of rows and columns and solving the PDE only a small number of times \cite{ballani2015hierarchical}.

Our central thesis is that cross approximation can be extended to matrices over a Hilbert space and used as a non-intrusive ROM method to approximate the solution map itself, rather than merely an associated scalar QoI.

\subsection{Contributions and outline}
We formally introduce our objects of study, \emph{Bochner matrices}, in \Cref{sec:bmx} and then analyze them through four distinct algebraic lenses: as linear operators (\Cref{sec:bmx_op}), as quasimatrices (\Cref{sec:bmx_qm}), as matrices (\Cref{sec:bmx_m}), and as third-order tensors (\Cref{sec:bmx_3t}).
Following this progression, we identify how the algebraic properties of Bochner matrices emerge from their structure, providing the necessary foundation for cross approximation.

We investigate the theoretical properties of Bochner cross approximation in \Cref{sec:cross}, deriving approximation guarantees based on (i)~the operator norm of the pseudoinverse of submatrices of singular factors and (ii)~the maximum-volume property.
Notably, we prove that the maximum-volume error bounds are smaller for Bochner matrices than for classical matrices of numbers (\Cref{subsec:rank_excess}).

\Cref{sec:numerics} presents an extension of \emph{adaptive cross approximation} \cite{bebendorf2000approximation, bebendorf2003adaptive} to Bochner matrices as a greedy algorithm for index selection.
We validate our algorithm, and the framework of Bochner cross approximation as a whole, as a non-intrusive ROM method in numerical experiments with parametric nonlinear Stokes equations.

Some concluding remarks are collected in \Cref{sec:remarks}.

\subsection{Notation}
We write vectors as $\Matrix{a,b}$ and matrices as $\Matrix{A,B}$, denoting the $n \times n$ identity matrix by $\Matrix{I}_n$.
The transpose is written as $\Matrix{A}^\trans$, the complex conjugate as $\overline{\Matrix{A}}$, and the conjugate transpose as $\Matrix{A}^\ast$.
We let $[n] = \{ 1, \ldots, n \}$.
By convention, the inner product $\Dotp{\cdot}{\cdot}{\Hilb}$ of a Hilbert space $\Hilb$ is linear in the second argument.
\section{Bochner matrices}
\label{sec:bmx}
Consider a Hilbert space $\Hilb$ over a field $\F \in \{ \Real, \Clx \}$.
Let $d \in \N$ and $n_1, \ldots, n_d \in \N$.
We say that an $n_1 \times \cdots \times n_d$ array with entries in $\Hilb$ is a \emph{Bochner tensor} of order $d$, which we refer to as a Bochner matrix when $d = 2$ and a Bochner vector when $d = 1$.
Equivalently, a Bochner tensor is a map
\begin{equation*}
    \Bensor{A} : [n_1] \times \cdots \times [n_d] \to \Hilb, \quad \mIndex{i} = (i_1, \ldots, i_d) \mapsto \Bensor{A}(i_1, \ldots, i_d) = \Bentry{a}_{\mIndex{i}}.
\end{equation*}
Tensors over the field $\F$ will be referred to simply as ``tensors,'' or as ``classical tensors'' for additional contrast.

Under entrywise addition and scalar multiplication, the set $\Hilb^{n_1 \times \cdots \times n_d}$ of Bochner tensors is a vector space.
We endow it with the family of $\ell_p(\Hilb)$ norms:
\begin{equation*}
    \Norm{\Bensor{A}}{\ell_p(\Hilb)} =
        \begin{cases}
            \left( \sum_{\mIndex{i} \in [n_1] \times \cdots \times [n_d]} \Norm{\Bentry{a}_{\mIndex{i}}}{\Hilb}^p \right)^{1/p}, & 1 \leq p < \infty, \\
            \max_{\mIndex{i} \in [n_1] \times \cdots \times [n_d]} \Norm{\Bentry{a}_{\mIndex{i}}}{\Hilb}, & p = \infty.
        \end{cases}
\end{equation*}
These norms are pairwise equivalent, satisfying the bounds
\begin{equation*}
    \Norm{\Bensor{A}}{\ell_q(\Hilb)} \leq \Norm{\Bensor{A}}{\ell_{p}(\Hilb)} \leq (n_1 \cdots n_d)^{\frac{1}{p} - \frac{1}{q}} \Norm{\Bensor{A}}{\ell_q(\Hilb)}, \quad p \leq q,
\end{equation*}
and satisfy H\"older's inequality:
\begin{equation*}
    |\Dotp{\Bensor{A}}{\Bensor{B}}{\ell_2(\Hilb)}| \leq \Norm{\Bensor{A}}{\ell_p(\Hilb)} \Norm{\Bensor{B}}{\ell_{p_\ast}(\Hilb)}, \quad p_{\ast} = \tfrac{p}{p-1},
\end{equation*}
where $\Dotp{\Bensor{A}}{\Bensor{B}}{\ell_2(\Hilb)} = \sum_{\mIndex{i}} \Dotp{\Bentry{a}_{\mIndex{i}}}{\Bentry{b}_{\mIndex{i}}}{\Hilb}$ is the inner product inducing the $\ell_2(\Hilb)$ norm.

A basic Cauchy-sequence argument shows that the vector space of Bochner tensors equipped with the $\ell_p(\Hilb)$ norm is complete, and so is itself a Hilbert space when $p = 2$.
Let us remark in passing that Bochner tensors can be identified with functions from Bochner--Lebesgue spaces with respect to the discrete counting measure.
\section{Bochner matrices as linear operators}
\label{sec:bmx_op}
An $m \times n$ Bochner matrix $\Batrix{A}$ can be multiplied from the right by classical vectors in $\F^n$.
This induces a linear operator $\Matrix{x} \mapsto \Batrix{A} \Matrix{x}$ between Hilbert spaces $\F^{n}$ and $\Hilb^{m}$.
The operator norm of $\Batrix{A}$ is equivalent to the $\ell_p(\Hilb)$ norms, satisfying the bounds $\Norm{\Batrix{A}}{\ell_{\infty}(\Hilb)} \leq \Norm{\Batrix{A}}{\ell_2 \to \ell_2(\Hilb)} \leq \Norm{\Batrix{A}}{\ell_1(\Hilb)}$, and thus the operator is bounded.
Furthermore, it has a finite rank not exceeding $n$, which we denote by $\rank{\Batrix{A}}$ and that satisfies the following basic property.

\begin{lemma}
\label{lemma:rank_AB}
Let $\Batrix{A} \in \Hilb^{m \times n}$ and $\Matrix{B} \in \F^{n \times k}$.
Then
\begin{equation*}
    \rank{\Batrix{A} \Matrix{B}} \leq \min\{ \rank{\Batrix{A}}, \rank{\Matrix{B}} \}, \quad \rank{\Batrix{A} \Matrix{B}} = \begin{cases}
            \rank{\Batrix{A}}, & \rank{\Matrix{B}} = n, \\
            \rank{\Matrix{B}}, & \rank{\Batrix{A}} = n.
        \end{cases}
\end{equation*}
\end{lemma}
\begin{proof}
    The inequality is standard. 
    The equalities follow directly from Sylvester's rank inequality, which applies as the operator of $\Batrix{A}$ has a finite-dimensional domain.
\end{proof}

As a bounded operator, $\Batrix{A}$ admits a bounded \emph{adjoint} $\Batrix{A}^{\ast} : \Hilb^{m} \to \F^{n}$ given by
\begin{equation*}
    \Batrix{A}^\ast \Batrix{y} = \begin{bmatrix}
        \Dotp{\Batrix{A}(:,1)}{\Batrix{y}}{\ell_2(\Hilb)} &
        \cdots &
        \Dotp{\Batrix{A}(:,n)}{\Batrix{y}}{\ell_2(\Hilb)}
    \end{bmatrix}^\trans.
\end{equation*}
We say that $\Batrix{A}$ has \emph{orthonormal columns} when its induced operator is a partial isometry, i.e., $\Batrix{A}^{\ast} \Batrix{A} = \Matrix{I}_n$, where we adopt the convention that the adjoint $\Batrix{A}^{\ast}$ can be applied to Bochner matrices columnwise.
Since the range of the operator induced by $\Batrix{A}$ is finite-dimensional (and hence closed), there exists a unique bounded \emph{pseudoinverse} $\Batrix{A}^{\pinv} : \Hilb^{m} \to \F^{n}$ defined by the standard Moore--Penrose equations \cite[Ch.~9]{ben2003generalized}.

Since $\Batrix{A}$ corresponds to a finite-rank (and hence compact) operator, it admits a singular value decomposition (SVD) with $\rank{\Batrix{A}}$ uniquely defined nonzero singular values \cite[Ch.~2]{gohberg1969introduction}.
We formulate this classical result in terms of Bochner matrices.

\begin{theorem}
\label{theorem:svd}
Let $\Batrix{A} \in \Hilb^{m \times n}$ be nonzero and $r = \rank{\Batrix{A}}$.
There exists a unique $\Matrix{\Sigma} = \mathrm{diag}(\sigma_1, \ldots, \sigma_r) \in \Real^{r \times r}$ with $\sigma_1 \geq \cdots \geq \sigma_r > 0$, and there exist $\Batrix{U} \in \Hilb^{m \times r}$ and $\Matrix{V} \in \F^{n \times r}$ with orthonormal columns such that $\Batrix{A} = \Batrix{U} \Matrix{\Sigma} \Matrix{V}^\ast$.
Furthermore, the pseudoinverse is given by $\Batrix{A}^{\pinv} = \Matrix{V} \Matrix{\Sigma}^{-1} \Batrix{U}^{\ast}$.
\end{theorem}

The optimality of the truncated SVD for low-rank approximation of compact operators is well-established in the context of symmetrically normed operator ideals \cite[Ch.~3, Lem.~6.1]{gohberg1969introduction}. 
Finite-rank operators belong to all symmetrically normed operator ideals, and hence the Eckart--Young--Mirsky theorem holds for Bochner matrices.
A more direct proof of this result can be obtained via a reduction to classical matrices.
In the formulation of the theorem below, we set $\sigma_j(\Batrix{A}) = 0$ for all $j > \rank{\Batrix{A}}$.

\begin{theorem}
\label{theorem:mirsky}
Let $\Batrix{A} \in \Hilb^{m \times n}$.
Let $\phi$ be a symmetric gauge function\footnote{To simplify submultiplicative norm bounds, we assume throughout the text our symmetric gauge functions to be normalized as $\phi(1, 0, \ldots, 0) = 1$. For instance, such normalization is used in \cite{gohberg1969introduction}.} inducing a unitarily invariant norm $\Norm{\Batrix{A}}{\phi} = \phi(\sigma_1(\Batrix{A}), \ldots, \sigma_n(\Batrix{A}))$.
For every $1 \leq k < \rank{\Batrix{A}}$, the best approximation of $\Batrix{A}$ by a Bochner matrix of rank not exceeding $k$ is achieved by the $k$-truncated SVD, denoted $\truncate{\Batrix{A}}{k} = \Batrix{U}_k \Matrix{\Sigma}_k \Matrix{V}_k^\ast$:
\begin{equation*}
    \Norm{\Batrix{A} - \truncate{\Batrix{A}}{k}}{\phi} = \min \set{\Norm{\Batrix{A} - \Batrix{B}}{\phi}}{\rank{\Batrix{B}} \leq k}.
\end{equation*}
\end{theorem}

Clearly, the $\ell_2(\Hilb)$ norm and the operator $\ell_2 \to \ell_2(\Hilb)$ norms are unitarily invariant.
Let us also note that every such norm $\Norm{\cdot}{\phi}$ is equivalent to the $\ell_p(\Hilb)$ norms, a property which we state without proof for brevity.
\section{Bochner matrices as quasimatrices}
\label{sec:bmx_qm}
A \emph{quasimatrix} \cite{trefethen2010householder, townsend2015continuous, shustin2022semi} is a generalized matrix with columns that are functions in $\Ell^2(0,1)$.
Adopting a broader view, the columns can lie in an abstract Hilbert space, which allows us to treat a Bochner matrix $\Batrix{A} \in \Hilb^{m \times n}$ as a quasimatrix with columns in $\Hilb^{m}$.
Consequently, $\Batrix{A}$ has exactly $\rank{\Batrix{A}}$ linearly independent columns and admits an \emph{interpolative decomposition}.

\begin{proposition}
\label{proposition:1s_interpol}
Let $\Batrix{A} \in \Hilb^{m \times n}$ be nonzero and $r = \rank{\Batrix{A}}$.
There exists an index set $\Index{J} \subseteq [n]$ of cardinality $r$ together with $\Matrix{B} \in \F^{r \times n}$ such that $\Batrix{A} = \Batrix{A}(:, \Index{J}) \Matrix{B}$.
\end{proposition}

As the columns of $\Batrix{A}$ belong to a Hilbert space, they can be orthogonalized via the Gram--Schmidt procedure, yielding a \emph{thin} QR decomposition.
The QR decomposition also provides a practical way to compute the SVD of a Bochner matrix, since it reduces the original problem to computing the SVD of a classical triangular matrix; the $n \times n$ Gram matrix of $\Batrix{A}$ could also be used to approach the computation of the SVD.

\begin{theorem}
\label{theorem:qr}
Let $\Batrix{A} \in \Hilb^{m \times n}$ satisfy $\rank{\Batrix{A}} = n$.
Then there exists a unique $\Batrix{Q} \in \Hilb^{m \times n}$ with orthonormal columns, and a unique upper triangular $\Matrix{R} \in \F^{n \times n}$ with strictly positive diagonal entries, such that $\Batrix{A} = \Batrix{Q} \Matrix{R}$.
\end{theorem}

\begin{corollary}
\label{corollary:qr_pivoted}
Let $\Batrix{A} \neq 0$ and $r = \rank{\Batrix{A}} < n$. 
There exist $\Batrix{Q} \in \Hilb^{m \times r}$ with orthonormal columns, an upper triangular $\Matrix{R}_1 \in \F^{r \times r}$ with strictly positive diagonal entries, $\Matrix{R}_2 \in \F^{r \times (n-r)}$, and an $n \times n$ permutation matrix $\Matrix{\Pi}$ such that $\Batrix{A}\Matrix{\Pi} = \Batrix{Q} [ \Matrix{R}_1~\Matrix{R}_2 ]$.
\end{corollary}
\begin{proof}
    Apply \cref{theorem:qr} to $\Batrix{A}(:, \Index{J})$ from \cref{proposition:1s_interpol}.
\end{proof}

The columnwise structure of Bochner matrices permits their low-rank approximation via \emph{approximate interpolative decompositions} as in \cref{proposition:1s_interpol}.
Given selected column indices $\Index{J}$, the approximation can be obtained as $\Batrix{A} \approx \Batrix{A}(:, \Index{J}) \Batrix{A}(:, \Index{J})^{\pinv} \Batrix{A}$.

\begin{proposition}
\label{proposition:cca_exact}
Let $\Batrix{A} \in \Hilb^{m \times n}$ and $\Index{J} \subseteq [n]$, and denote $\Batrix{C} = \Batrix{A}(:, \Index{J})$.
Then $\Batrix{A} = \Batrix{C} \Batrix{C}^{\pinv} \Batrix{A}$ if and only if $\rank{\Batrix{C}} = \rank{\Batrix{A}}$.
\end{proposition}
\begin{proof}
Follows from \cref{lemma:rank_AB} and the fact that $\Batrix{C} \Batrix{C}^{\pinv}$ is an orthogonal projection operator onto the column span of $\Batrix{C}$.
\end{proof}

When the goal is to select a small number of columns, the \emph{column subset selection} problem aims to identify the optimal ones and can be approached via rank-revealing QR decompositions or volume-based sampling \cite{kressner2025approximation}---relying only on the inner-product structure of the space of columns.

The \emph{rows} of a quasimatrix remain undefined unless its columns possess additional structure, e.g., its column-functions having continuous representatives in $\Ell^2(0,1)$ \cite{townsend2015continuous} or lying in a reproducing kernel Hilbert space \cite{shustin2022semi}.
Bochner matrices exhibit a different enabling structure.
\section{Bochner matrices as matrices}
\label{sec:bmx_m}
Finally, we acknowledge that Bochner matrices are two-dimensional arrays, and hence can be transposed.
However, the ranks of $\Batrix{A}$ and $\Batrix{A}^{\trans}$ are, in general, distinct.
An example demonstrating this divergence is an $m \times n$ Bochner matrix whose entries constitute an orthonormal system in $\Hilb$: it has $\rank{\Batrix{A}} = n$ and $\rank{\Batrix{A}^\trans} = m$.
Therefore, we distinguish between the \emph{column rank} and \emph{row rank} of a Bochner matrix:
\begin{equation*}
    \Rank{c}{\Batrix{A}} = \rank{\Batrix{A}}, \quad \Rank{r}{\Batrix{A}} = \rank{\Batrix{A}^\trans}.
\end{equation*}
All of the decompositions discussed so far apply equally well to transposed Bochner matrices, though it is important to keep in mind that the like decompositions of $\Batrix{A}$ and $\Batrix{A}^\trans$ are not connected in any immediate way.
Another property of classical matrices that fails for Bochner matrices is the commutativity of transposition with taking the adjoint or pseudoinverse, because the latter are defined strictly as linear operators.

Since we now allow ourselves to access the rows, we can multiply Bochner matrices on the left by classical matrices.
While \cref{lemma:rank_AB} describes the row rank of $\Matrix{B} \Batrix{A}$ via transposition, its column rank needs to be treated separately.

\begin{lemma}
\label{lemma:rank_BA}
Let $\Batrix{A} \in \Hilb^{m \times n}$ and $\Matrix{B} \in \F^{k \times m}$.
If $\rank{\Matrix{B}} = m$ then $\rank{\Matrix{B} \Batrix{A}} = \rank{\Batrix{A}}$, otherwise $\rank{\Matrix{B} \Batrix{A}} \leq \rank{\Batrix{A}}$.
\end{lemma}
\begin{proof}
Applying the rank-nullity theorem to the linear operator $\Matrix{B} \Batrix{A} : \F^n \to \Hilb^k$ and representing its kernel as a direct sum $\ker(\Matrix{B} \Batrix{A}) = \ker(\Batrix{A}) \oplus (\ker(\Batrix{A})^\perp \cap \ker(\Matrix{B} \Batrix{A}))$, we obtain the following identity:
\begin{equation*}
    \rank{\Matrix{B} \Batrix{A}} = \rank{\Batrix{A}} - \dim(\ker(\Batrix{A})^\perp \cap \ker(\Matrix{B} \Batrix{A})).
\end{equation*}
The rank inequality is its immediate consequence.
To prove the equality, it is sufficient to show that $\ker(\Batrix{A}) = \ker(\Matrix{B} \Batrix{A})$.
Let $\Batrix{z} \in \Hilb^m$ satisfy $\Matrix{B} \Batrix{z} = 0$.
For any $\Bentry{h} \in \Hilb$, we have
\begin{equation*}
    \begin{bmatrix}
        0 \\
        \vdots \\
        0
    \end{bmatrix} = 
    \begin{bmatrix}
        \Dotp{\sum_{j = 1}^{m} b_{1,j} \Bentry{z}_j}{\Bentry{h}}{\Hilb} \\
        \vdots \\
        \Dotp{\sum_{j = 1}^{m} b_{k,j} \Bentry{z}_j}{\Bentry{h}}{\Hilb}
    \end{bmatrix} = 
    \Matrix{B} \begin{bmatrix}
        \Dotp{\Bentry{z}_1}{\Bentry{h}}{\Hilb} \\
        \vdots \\
        \Dotp{\Bentry{z}_m}{\Bentry{h}}{\Hilb}
    \end{bmatrix}.
\end{equation*}
Since $\Matrix{B}$ has full column rank by assumption, $\Dotp{\Bentry{z}_j}{\Bentry{h}}{\Hilb} = 0$ for each $j$.
But as $\Bentry{h}$ is arbitrary, we conclude that $\Batrix{z} = 0$ and the kernels coincide.
\end{proof}

Furthermore, the access to rows provides us with a different way to construct the approximate interpolative decomposition, generalizing \cref{proposition:cca_exact}.

\begin{proposition}
\label{proposition:cgr_exact}
Let $\Batrix{A} \in \Hilb^{m \times n}$, $\Index{I} \subseteq [m]$, and $\Index{J} \subseteq [n]$.
Denote $\Batrix{C} = \Batrix{A}(:, \Index{J})$, $\Batrix{R} = \Batrix{A}(\Index{I}, :)$, and $\Batrix{G} = \Batrix{A}(\Index{I}, \Index{J})$.
Then $\Batrix{A} = \Batrix{C} \Batrix{G}^{\pinv} \Batrix{R}$ if and only if $\rank{\Batrix{G}} = \rank{\Batrix{A}}$.
\end{proposition}
\begin{proof}
If $\Batrix{A} = \Batrix{C} \Batrix{G}^\pinv \Batrix{R}$ then $\rank{\Batrix{A}} \leq \rank{\Batrix{G}^\pinv \Batrix{R}} \leq \rank{\Batrix{G}} \leq \rank{\Batrix{A}}$.
Conversely, we have $\rank{\Batrix{A}} = \rank{\Batrix{G}} \leq \rank{\Batrix{C}} \leq \rank{\Batrix{A}}$.
By \cref{proposition:cca_exact}, $\Batrix{A} = \Batrix{C} \Batrix{C}^\pinv \Batrix{A}$ and thus $\Batrix{R} = \Batrix{G} \Batrix{C}^\pinv \Batrix{A}$.
On the other hand, $\rank{\Batrix{A}} = \rank{\Batrix{G}} \leq \rank{\Batrix{R}} \leq \rank{\Batrix{A}}$, so that $\Batrix{R} = \Batrix{G} \Batrix{G}^\pinv \Batrix{R}$ by \cref{proposition:cca_exact}.
Let $\Batrix{C} = \Batrix{C}(:, \Index{K}) \Matrix{B}$ be an interpolative decomposition from \cref{proposition:1s_interpol}, and note that $\Matrix{B}$ has full row rank.
We then have $\Batrix{G} = \Batrix{G}(:, \Index{K}) \Matrix{B}$, and $\Batrix{G}(:, \Index{K})$ has full column rank by \cref{lemma:rank_AB}.
To finish the proof,
\begin{equation*}
    \Batrix{G}^\pinv \Batrix{R} = \Batrix{G}^\pinv \Batrix{G} \Batrix{C}^\pinv \Batrix{A} = \Matrix{B}^\pinv \Batrix{G}(:, \Index{K})^\pinv \Batrix{G}(:, \Index{K}) \Matrix{B} \Matrix{B}^\pinv \Batrix{C}(:, \Index{K})^\pinv \Batrix{A} = \Matrix{B}^\pinv \Batrix{C}(:, \Index{K})^\pinv \Batrix{A} = \Batrix{C}^\pinv \Batrix{A}.
\end{equation*}
\end{proof}

The interpolative decomposition in \cref{proposition:cgr_exact} is written in exactly the same form as the skeleton decomposition of classical matrices \cite{goreinov1997theory}, i.e., $\Matrix{A}(:, \Index{J}) \Matrix{A}(\Index{I}, \Index{J})^{\pinv} \Matrix{A}(\Index{I}, :)$.
Despite this syntactic similarity, there are two crucial differences.
First, the classical skeleton decomposition is a product of three factors, whereas the decomposition in \cref{proposition:cgr_exact} fundamentally collapses into two factors because the pseudoinverse operator cannot stand alone as a purely algebraic middle factor.
Second, the skeleton decomposition of classical matrices is a CUR decomposition \cite{hamm2020perspectives}; that is, every column is expressed as a linear combination of selected columns and every row is expressed as a linear combination of selected rows.
In contrast, the decomposition in \cref{proposition:cgr_exact} is interpolative only with respect to the columns and fails to leverage the row rank.

In fact, CUR decompositions in the sense described above are structurally impossible for Bochner matrices: they would contain two Bochner-matrix factors, yet those cannot be multiplied within our framework.
In this regard, consider a verbatim extension of the projection-based classical CUR decomposition $\Matrix{C} \Matrix{C}^{\pinv} \Matrix{A} \Matrix{R}^{\pinv} \Matrix{R}$ \cite{mahoney2009cur}, which we interpret not as a three-factor product with $\Matrix{C}^{\pinv} \Matrix{A} \Matrix{R}^{\pinv}$ in the middle, but as the simultaneous application of column and row orthogonal projectors.
For Bochner matrices, we can write this down as $(\Batrix{R}^\trans (\Batrix{R}^\trans)^\pinv (\Batrix{C} \Batrix{C}^\pinv \Batrix{A})^\trans)^\trans$ or $\Batrix{C} \Batrix{C}^\pinv (\Batrix{R}^\trans (\Batrix{R}^\trans)^\pinv \Batrix{A}^\trans)^\trans$.
Consider a simple $2 \times 2$ example where the entries of $\Batrix{A}$ form an orthonormal system in $\Hilb$, and let $\Index{I} = \Index{J} = \{ 1 \}$.
Then the two projection-based CUR extensions evaluate to
\begin{equation*}
    (\Batrix{R}^\trans (\Batrix{R}^\trans)^\pinv (\Batrix{C} \Batrix{C}^\pinv \Batrix{A})^\trans)^\trans = \frac{1}{2} \begin{bmatrix}
        \Batrix{A}(1,1) & \Batrix{A}(1,2) \\
        0 & 0
    \end{bmatrix}, \quad
    \Batrix{C} \Batrix{C}^\pinv (\Batrix{R}^\trans (\Batrix{R}^\trans)^\pinv \Batrix{A}^\trans)^\trans = \frac{1}{2} \begin{bmatrix}
        \Batrix{A}(1,1) & 0 \\
        \Batrix{A}(2,1) & 0
    \end{bmatrix}.
\end{equation*}
They do not agree, the column and row projections overriding each other, and neither is strictly a CUR decomposition.
Note also that the subspaces we selected to project on are dominant left singular subspaces of $\Batrix{A}$ and $\Batrix{A}^\trans$.
\section{Bochner matrices as third-order tensors}
\label{sec:bmx_3t}
As demonstrated, a CUR-style application of the one-sided interpolative decompositions from \cref{proposition:cca_exact,proposition:cgr_exact} does not yield a simultaneous reduction of both the column and row ranks.
Instead, we will extend classical \emph{two-sided interpolative decompositions} \cite{voronin2017efficient} into the Bochner setting.

\begin{proposition}
\label{proposition:2s_interpol}
Let $\Batrix{A} \in \Hilb^{m \times n}$ be nonzero with $\rho = \Rank{r}{\Batrix{A}}$, $\kappa = \Rank{c}{\Batrix{A}}$.
There exist index sets $\Index{I} \subseteq [m]$ of cardinality $\rho$ and $\Index{J} \subseteq [n]$ of cardinality $\kappa$, together with $\Matrix{B} \in \F^{m \times \rho}$ and $\Matrix{C} \in \F^{\kappa \times n}$, such that $\Batrix{A} = \Matrix{B} \Batrix{A}(\Index{I}, \Index{J}) \Matrix{C}$.
\end{proposition}
\begin{proof}
Applying \cref{proposition:1s_interpol} to $\Batrix{A}$ and $\Batrix{A}^\trans$, we get $\Batrix{A} = \Batrix{A}(:, \Index{J}) \Matrix{C} = \Matrix{B} \Batrix{A}(\Index{I}, :)$.
\end{proof}

The structure of the decomposition in \cref{proposition:2s_interpol}, having a Bochner matrix as its inner core and classical matrices as its outer factors, strictly mirrors the Tucker decomposition of classical tensors \cite{ballard2025tensor}.
Indeed, in the finite-dimensional special case of $\Hilb = \F^{k}$ this reduces to the classical Tucker2 decomposition; and in the general case, the space of Bochner matrices $\Hilb^{m \times n}$ is isometrically isomorphic to the tensor-product Hilbert space $\Hilb \krp \F^{m} \krp \F^{n}$ \cite{hackbusch2019tensor}.
Motivated by these correspondences, we call every decomposition of the form $\Batrix{A} = \Matrix{B} \Batrix{G} \Matrix{D}$ a \emph{Tucker decomposition} and refer to it as a \emph{cross decomposition} when its core $\Batrix{G}$ is a submatrix of $\Batrix{A}$.
Accordingly, we define the \emph{Tucker rank} of a Bochner matrix as a tuple $\Rank{T}{\Batrix{A}} = (\Rank{r}{\Batrix{A}}, \Rank{c}{\Batrix{A}})$.

The following theorem extends \cref{proposition:cgr_exact} to cross decompositions and connects them to the CUR-style application of one-sided interpolative decompositions.
This generalizes properties of classical matrices and tensors derived in \cite{caiafa2010generalizing, hamm2020perspectives, cai2021mode} to Bochner matrices.

\begin{theorem}
\label{theorem:cross_exact}
Let $\Batrix{A} \in \Hilb^{m \times n}$, $\Index{I} \subseteq [m]$, and $\Index{J} \subseteq [n]$.
Denote $\Batrix{C} = \Batrix{A}(:, \Index{J})$, $\Batrix{R} = \Batrix{A}(\Index{I},:)$, and $\Batrix{G} = \Batrix{A}(\Index{I}, \Index{J})$. Then the following are equivalent:
\begin{enumerate}
    \item $\Rank{c}{\Batrix{C}} = \Rank{c}{\Batrix{A}}$ and $\Rank{r}{\Batrix{R}} = \Rank{r}{\Batrix{A}}$;
    \item $\Rank{T}{\Batrix{G}} = \Rank{T}{\Batrix{A}}$;
    \item $\Batrix{A} = (\Batrix{R}^\trans (\Batrix{R}^\trans)^\pinv (\Batrix{C} \Batrix{C}^\pinv \Batrix{A})^\trans)^\trans = \Batrix{C} \Batrix{C}^\pinv (\Batrix{R}^\trans (\Batrix{R}^\trans)^\pinv \Batrix{A}^\trans)^\trans$;
    \item $\Batrix{A} = ((\Batrix{G}^\trans)^\pinv \Batrix{C}^\trans)^\trans \Batrix{G} (\Batrix{G}^\pinv \Batrix{R})$.
\end{enumerate}
\end{theorem}
\begin{proof}
$\mathbf{(1 \Rightarrow 2).}$
\cref{proposition:cca_exact} gives $\Batrix{A} = \Batrix{C} \Batrix{C}^\pinv \Batrix{A}$, and thus $\Batrix{R} = \Batrix{G} \Batrix{C}^\pinv \Batrix{A}$.
It follows that $\Rank{c}{\Batrix{G}} = \Rank{c}{\Batrix{R}}$ since $\Batrix{G}$ is a submatrix of $\Batrix{R}$.
Similarly, $\Rank{r}{\Batrix{G}} = \Rank{r}{\Batrix{C}}$, and there exists $\Index{K} \subseteq [m]$ together with $\Matrix{X}$ of full column rank such that $\Batrix{C} = \Matrix{X} \Batrix{G}(\Index{K}, :)$.
By \cref{lemma:rank_BA}, we get $\Rank{c}{\Batrix{C}} = \Rank{c}{\Batrix{G}(\Index{K}, :)} \leq \Rank{c}{\Batrix{G}} \leq \Rank{c}{\Batrix{A}} = \Rank{c}{\Batrix{C}}$.
By analogy, $\Rank{r}{\Batrix{G}} = \Rank{r}{\Batrix{A}}$.

$\mathbf{(2 \Rightarrow 1).}$
We have $\Rank{c}{\Batrix{A}} = \Rank{c}{\Batrix{G}} \leq \Rank{c}{\Batrix{C}} \leq \Rank{c}{\Batrix{A}}$ as $\Batrix{G}$ is a submatrix of $\Batrix{C}$.
Similarly, $\Rank{r}{\Batrix{R}} \leq \Rank{r}{\Batrix{A}}$.

$\mathbf{(1 \Rightarrow 3).}$
By \cref{proposition:cca_exact}, $\Batrix{A} = \Batrix{C} \Batrix{C}^\pinv \Batrix{A} = (\Batrix{R}^\trans (\Batrix{R}^\trans)^\pinv \Batrix{A}^\trans)^\trans = \Batrix{C} \Batrix{C}^\pinv (\Batrix{R}^\trans (\Batrix{R}^\trans)^\pinv \Batrix{A}^\trans)^\trans$.
Repeat the argument for $\Batrix{A}^\trans$ to get the second representation.

$\mathbf{(3 \Rightarrow 1).}$
By \cref{lemma:rank_AB}, $\Rank{c}{\Batrix{A}} \leq \Rank{c}{\Batrix{C}} \leq \Rank{c}{\Batrix{A}}$, and similarly for $\Batrix{R}^\trans$. 

$\mathbf{(2 \Rightarrow 4).}$
\cref{proposition:cgr_exact} yields $\Batrix{A} = \Batrix{C} \Batrix{G}^\pinv \Batrix{R} = (\Batrix{R}^\trans (\Batrix{G}^\trans)^\pinv \Batrix{C}^\trans)^\trans$, whence we obtain $\Batrix{C}^\trans = \Batrix{G}^\trans (\Batrix{G}^\trans)^\pinv \Batrix{C}^\trans$.

$\mathbf{(4 \Rightarrow 2).}$
Note that $\Batrix{C} = ((\Batrix{G}^\trans)^\pinv \Batrix{C}^\trans)^\trans \Batrix{G} (\Batrix{G}^\pinv \Batrix{G}) = ((\Batrix{G}^\trans)^\pinv \Batrix{C}^\trans)^\trans \Batrix{G}$ by the definition of pseudoinverse.
Then $\Batrix{A} = \Batrix{C} \Batrix{G}^\pinv \Batrix{R}$, and \cref{proposition:cgr_exact} gives $\Rank{c}{\Batrix{G}} = \Rank{c}{\Batrix{A}}$.
The same argument gives $\Rank{r}{\Batrix{G}} = \Rank{r}{\Batrix{A}}$ when applied to $\Batrix{A}^\trans$.
\end{proof}

The interpretation of Bochner matrices as third-order tensors paves the way for their low-rank approximation by means of \emph{higher-order SVD} (HOSVD) \cite{de2000multilinear}.
However, as the example at the end of \Cref{sec:bmx_m} shows, unlike the case of classical tensors, low-rank truncation for each unfolding (i.e., $\Batrix{A}$ and $\Batrix{A}^\trans$) cannot be performed via column-space projection: the two projections do not commute and interfere with one another.
Instead, we must apply classical orthogonal projection matrices on the right, thereby projecting in the row space.
Such construction is perfectly aligned with the definition of mode unfoldings in abstract tensor-product Hilbert spaces \cite{hackbusch2019tensor}.

\begin{theorem}
\label{theorem:hosvd}
Let $\Batrix{A} \in \Hilb^{m \times n}$ be nonzero, and let $\Batrix{A} = \Batrix{U} \Matrix{\Sigma} \Matrix{V}^\ast$ and $\Batrix{A}^\trans = \Batrix{P} \Matrix{\Lambda} \Matrix{Q}^\ast$ be its SVDs.
For every pair of $1 \leq \rho \leq \Rank{r}{\Batrix{A}}$ and $1 \leq \kappa \leq \Rank{c}{\Batrix{A}}$, denote by $\truncate{\Batrix{A}}{\kappa} = \Batrix{U}_{\kappa} \Matrix{\Sigma}_{\kappa} \Matrix{V}_{\kappa}^\ast$ and $\truncate{\Batrix{A}^\trans}{\rho} = \Batrix{P}_{\rho} \Matrix{\Lambda}_{\rho} \Matrix{Q}_{\rho}^\ast$ the corresponding truncated SVDs.
Then
\begin{align*}
    \Norm{\Batrix{A} - \overline{\Matrix{Q}}_\rho \Matrix{Q}_\rho^\trans \Batrix{A} \Matrix{V}_\kappa \Matrix{V}_\kappa^\ast}{\ell_2(\Hilb)} &\leq \sqrt{\Norm{\Matrix{\Sigma} - \Matrix{\Sigma}_\kappa}{\ell_2}^2 + \Norm{\Matrix{\Lambda} - \Matrix{\Lambda}_\rho}{\ell_2}^2} \\
    &\leq \sqrt{2} \min \set[\Big]{\Norm{\Batrix{A} - \Batrix{B}}{\ell_2(\Hilb)}}{\Rank{r}{\Batrix{B}} \leq \rho,~\Rank{c}{\Batrix{B}} \leq \kappa}.
\end{align*}
\end{theorem}
\begin{proof}
The Pythagorean theorem gives
\begin{align*}
    \Norm{\Batrix{A} - \overline{\Matrix{Q}}_\rho \Matrix{Q}_\rho^\trans \Batrix{A} \Matrix{V}_\kappa \Matrix{V}_\kappa^\ast}{\ell_2(\Hilb)}^2 &= \Norm{\Batrix{A} (\Matrix{I} - \Matrix{V}_\kappa \Matrix{V}_\kappa^\ast)}{\ell_2(\Hilb)}^2 + \Norm{(\Matrix{I} - \overline{\Matrix{Q}}_\rho \Matrix{Q}_\rho^\trans) \Batrix{A} \Matrix{V}_\kappa \Matrix{V}_\kappa^\ast}{\ell_2(\Hilb)}^2 \\
    &\leq \Norm{\Batrix{A} (\Matrix{I} - \Matrix{V}_\kappa \Matrix{V}_\kappa^\ast)}{\ell_2(\Hilb)}^2 + \Norm{\Batrix{A}^\trans (\Matrix{I} - \Matrix{Q}_\rho \Matrix{Q}_\rho^\ast)}{\ell_2(\Hilb)}^2 \\
    &= \Norm{\Matrix{\Sigma} - \Matrix{\Sigma}_\kappa}{\ell_2}^2 + \Norm{\Matrix{\Lambda} - \Matrix{\Lambda}_\rho}{\ell_2}^2.
\end{align*}
A standard corollary of the Eckart--Young--Mirsky theorem (\cref{theorem:mirsky}) is that the set of Bochner matrices with bounded Tucker rank is closed, and therefore the minimal distance is achieved at some $\Batrix{B}_{\circ}$.
Invoking \cref{theorem:mirsky} again, we obtain
\begin{equation*}
    \Norm{\Matrix{\Sigma} - \Matrix{\Sigma}_\kappa}{\ell_2}^2 \leq \Norm{\Batrix{A} - \Batrix{B}_{\circ}}{\ell_2(\Hilb)}^2, \quad \Norm{\Matrix{\Lambda} - \Matrix{\Lambda}_\rho}{\ell_2}^2 \leq \Norm{\Batrix{A}^\trans - \Batrix{B}_{\circ}^\trans}{\ell_2(\Hilb)}^2.
\end{equation*}
\end{proof}

The proof is a direct adaptation of the HOSVD argument for classical tensors \cite{de2000multilinear}, and it extends to tensors in abstract tensor-product Hilbert spaces because their mode unfoldings are Hilbert--Schmidt operators (even though \cite[Thm.~10.3]{hackbusch2019tensor} is restricted to classical tensors). 
A similar quasioptimality bound can be shown for the \emph{sequentially truncated HOSVD} of Bochner matrices (see \cite{vannieuwenhoven2012new} and \cite[Thm.~10.5]{hackbusch2019tensor}).
\section{Cross approximation of Bochner matrices}
\label{sec:cross}
In the preceding sections, we established that Bochner matrices can be approximated with low Tucker rank via the truncated HOSVD (\cref{theorem:hosvd}) and derived necessary and sufficient conditions for their cross decompositions to be exact (\cref{theorem:cross_exact}).
However, our motivating application introduced in \cref{sec:intro:rom} relies on \emph{approximate} cross decompositions, and therefore this section investigates the theoretical properties of \emph{cross approximation}:
\begin{equation*}
    \cross{\Batrix{B}}{\Index{I}, \Index{J}} = ((\Batrix{G}_{\Batrix{B}}^\trans)^\pinv \Batrix{C}_{\Batrix{B}}^\trans)^\trans \Batrix{G}_{\Batrix{B}} (\Batrix{G}_{\Batrix{B}}^\pinv \Batrix{R}_{\Batrix{B}}),
\end{equation*}
where $\Batrix{C}_{\Batrix{B}} = \Batrix{B}(:, \Index{J})$, $\Batrix{G}_{\Batrix{B}} = \Batrix{B}(\Index{I}, \Index{J})$, $\Batrix{R}_{\Batrix{B}} = \Batrix{B}(\Index{I}, :)$ are the selected submatrices.
We will drop the subscripts for the submatrices of $\Batrix{A}$ and use regular font for the submatrices of classical matrices.

\subsection{Interpolation property}
When the conditions of \cref{theorem:cross_exact} are not met and cross approximation fails to recover the entire Bochner matrix, it still interpolates entries residing in the selected rows and columns.
The following theorem formally extends the interpolation property of classical cross approximation to Bochner matrices.
For instance, it shows that every cross approximant is a cross decomposition of itself.

\begin{theorem}
\label{theorem:cross_interpol}
Let $\Batrix{A} \in \Hilb^{m \times n}$ and $\Batrix{X} = \cross{\Batrix{A}}{\Index{I}, \Index{J}}$.
Then
\begin{enumerate}
    \item $\Batrix{X}(\Index{I}, \Index{J}) = \Batrix{G}$ and $\Rank{T}{\Batrix{X}} = \Rank{T}{\Batrix{G}}$;
    \item $\Batrix{X}(\Index{I},:) = \Batrix{R}$ if and only if $\Rank{c}{\Batrix{G}} = \Rank{c}{\Batrix{R}}$;
    \item $\Batrix{X}(:,\Index{J}) = \Batrix{C}$ if and only if $\Rank{r}{\Batrix{G}} = \Rank{r}{\Batrix{C}}$;
    \item $\Batrix{X} = \Batrix{A}$ if and only if $\Rank{T}{\Batrix{G}} = \Rank{T}{\Batrix{A}}$.
\end{enumerate}
\end{theorem}
\begin{proof}
$\mathbf{1.}$
By the defining Moore--Penrose equations of pseudoinverse:
\begin{equation*}
    \Batrix{X}(\Index{I}, \Index{J}) = ((\Batrix{G}^\trans)^\pinv \Batrix{G}^\trans)^\trans \Batrix{G} (\Batrix{G}^\pinv \Batrix{G}) = ((\Batrix{G}^\trans)^\pinv \Batrix{G}^\trans)^\trans \Batrix{G} = (\Batrix{G}^\trans (\Batrix{G}^\trans)^\pinv \Batrix{G}^\trans)^\trans = \Batrix{G}.
\end{equation*}
By \cref{lemma:rank_AB,lemma:rank_BA}, $\Rank{T}{\Batrix{X}}$ is upper bounded by $\Rank{T}{\Batrix{G}}$.
As $\Batrix{G}$ is a submatrix of $\Batrix{X}$, the reverse inequality holds as well.

$\mathbf{2.}$
By definition, $\Batrix{X}(\Index{I},:) = \Batrix{G} \Batrix{G}^\pinv \Batrix{R}$.
\cref{proposition:cca_exact} guarantees that $\Batrix{X}(\Index{I},:) = \Batrix{R}$ if and only if $\Rank{c}{\Batrix{G}} = \Rank{c}{\Batrix{R}}$.

$\mathbf{3.}$
Apply the same argument to $\Batrix{X}^\trans$.

$\mathbf{4.}$ Follows from \cref{theorem:cross_exact}.
\end{proof}

\subsection{Error bound for one-sided interpolation}
As an intermediate step towards the error bound for Bochner cross approximation, we derive an error bound for the approximate one-sided interpolation (\cref{proposition:cgr_exact}).
The proof mirrors the argument used for cross approximation of classical perturbed low-rank matrices \cite{hamm2021perturbations}.
We begin with a technical lemma.

\begin{lemma}
\label{lemma:sub_pinv}
Let $\Batrix{B} \in \Hilb^{m \times n}$ be nonzero and satisfy $\Rank{c}{\Batrix{G}_{\Batrix{B}}} = \Rank{c}{\Batrix{B}}$, and let $\Batrix{B} = \Batrix{U} \Matrix{\Sigma} \Matrix{V}^\ast$ be its SVD.
Then
\begin{equation*}
    \Norm{\Batrix{C}_{\Batrix{B}} \Batrix{G}_{\Batrix{B}}^\pinv}{\ell_2(\Hilb) \to \ell_2(\Hilb)} = \Norm{\Batrix{R}_{\Batrix{U}}^\pinv}{\ell_2(\Hilb) \to \ell_2}, \quad \Norm{\Batrix{G}_{\Batrix{B}}^\pinv \Batrix{R}_{\Batrix{B}}}{\ell_2 \to \ell_2} = \Norm{\Matrix{C}_{\Matrix{V}^\ast}^\pinv}{\ell_2 \to \ell_2}.
\end{equation*}
\end{lemma}
\begin{proof}
By definition, $\Batrix{G}_{\Batrix{B}} = \Batrix{U}(\Index{I}, :) \Matrix{\Sigma} \Matrix{V}(\Index{J}, :)^\ast = \Batrix{R}_{\Batrix{U}} \Matrix{\Sigma} \Matrix{C}_{\Matrix{V}^\ast}$; likewise, $\Batrix{C}_{\Batrix{B}} = \Batrix{U} \Matrix{\Sigma} \Matrix{C}_{\Matrix{V}^\ast}$ and $\Batrix{R}_{\Batrix{B}} = \Batrix{R}_{\Batrix{U}} \Matrix{\Sigma} \Matrix{V}^\ast$.
The rank assumption guarantees, via \cref{lemma:rank_AB}, that $\Batrix{R}_{\Batrix{U}}$ has full column rank and $\Matrix{C}_{\Matrix{V}^\ast}$ has full row rank.
Therefore, $\Batrix{G}_{\Batrix{B}}^\pinv = \Matrix{C}_{\Matrix{V}^\ast}^\pinv \Matrix{\Sigma}^{-1} \Batrix{R}_{\Batrix{U}}^\pinv$ and thus
\begin{equation*}
    \Batrix{C}_{\Batrix{B}} \Batrix{G}_{\Batrix{B}}^\pinv = \Batrix{U} \Matrix{\Sigma} \Matrix{C}_{\Matrix{V}^\ast} \Matrix{C}_{\Matrix{V}^\ast}^\pinv \Matrix{\Sigma}^{-1} \Batrix{R}_{\Batrix{U}}^\pinv = \Batrix{U} \Batrix{R}_{\Batrix{U}}^\pinv, \quad \Batrix{G}_{\Batrix{B}}^\pinv \Batrix{R}_{\Batrix{B}} = \Matrix{C}_{\Matrix{V}^\ast}^\pinv \Matrix{\Sigma}^{-1} \Batrix{R}_{\Batrix{U}}^\pinv \Batrix{R}_{\Batrix{U}} \Matrix{\Sigma} \Matrix{V}^\ast = \Matrix{C}_{\Matrix{V}^\ast}^\pinv \Matrix{V}^\ast.
\end{equation*}
It remains to note that $\Batrix{U}$ and $\Matrix{V}^\ast$, being partial isometries, do not alter the respective operator norms.
\end{proof}

\begin{theorem}
\label{theorem:cgr_error}
Let $\Batrix{A} = \Batrix{B} + \Batrix{E} \in \Hilb^{m \times n}$ with $\Batrix{B} \neq 0$.
Let $\Batrix{B} = \Batrix{U} \Matrix{\Sigma} \Matrix{V}^\ast$ be an SVD, and let $\Rank{c}{\Batrix{G}_{\Batrix{B}}} = \Rank{c}{\Batrix{B}}$.
Let $\Norm{\cdot}{\phi}$ be a unitarily invariant norm induced by a symmetric gauge function $\phi$.
Denote $\eta_r = \Norm{\Batrix{R}_{\Batrix{U}}^\pinv}{\ell_2(\Hilb) \to \ell_2}$ and $\eta_c = \Norm{\Matrix{C}_{\Matrix{V}^\ast}^\pinv}{\ell_2 \to \ell_2}$.
Then
\begin{align*}
    \Norm{\Batrix{B} &- \Batrix{C} \Batrix{G}^\pinv \Batrix{R}}{\phi} \leq \eta_r \Norm{\Batrix{R}_{\Batrix{E}}}{\phi} + \eta_c \Norm{\Batrix{C}_{\Batrix{E}}}{\phi} + 3 \eta_r \eta_c \Norm{\Batrix{G}_{\Batrix{E}}}{\phi} \\
    &+ \Norm{\Batrix{G}^\pinv}{\ell_2(\Hilb) \to \ell_2} \Big( \Norm{\Batrix{R}_{\Batrix{E}}}{\phi} \Norm{\Batrix{C}_{\Batrix{E}}}{\phi} + \big(\eta_r \Norm{\Batrix{R}_{\Batrix{E}}}{\phi} + \eta_c \Norm{\Batrix{C}_{\Batrix{E}}}{\phi} + \eta_r \eta_c \Norm{\Batrix{G}_{\Batrix{E}}}{\phi}\big) \Norm{\Batrix{G}_{\Batrix{E}}}{\phi} \Big).
\end{align*}
\end{theorem}
\begin{proof}
The rank assumption together with \cref{proposition:cgr_exact} yields
\begin{equation*}
    \Norm{\Batrix{B} - \Batrix{C} \Batrix{G}^\pinv \Batrix{R}}{\phi} \leq \Norm{\Batrix{B} - \Batrix{C}_{\Batrix{B}} \Batrix{G}_{\Batrix{B}}^\pinv \Batrix{R}_{\Batrix{B}}}{\phi} + \Norm{\Batrix{C}_{\Batrix{B}} \Batrix{G}_{\Batrix{B}}^\pinv \Batrix{R}_{\Batrix{B}} - \Batrix{C} \Batrix{G}^\pinv \Batrix{R}}{\phi} = \Norm{\Batrix{C}_{\Batrix{B}} \Batrix{G}_{\Batrix{B}}^\pinv \Batrix{R}_{\Batrix{B}} - \Batrix{C} \Batrix{G}^\pinv \Batrix{R}}{\phi},
\end{equation*}
and therefore
\begin{align*}
    \Norm{\Batrix{B} - \Batrix{C} \Batrix{G}^\pinv \Batrix{R}}{\phi} &\leq \Norm{\Batrix{C}_{\Batrix{B}} \Batrix{G}_{\Batrix{B}}^\pinv \Batrix{R}_{\Batrix{B}} - \Batrix{C}_{\Batrix{B}} \Batrix{G}^\pinv \Batrix{R}_{\Batrix{B}}}{\phi} \\
    &+ \Norm{\Batrix{C}_{\Batrix{B}} \Batrix{G}^\pinv \Batrix{R}_{\Batrix{B}} - \Batrix{C}_{\Batrix{B}} \Batrix{G}^\pinv \Batrix{R}}{\phi} + \Norm{\Batrix{C}_{\Batrix{B}} \Batrix{G}^\pinv \Batrix{R} - \Batrix{C} \Batrix{G}^\pinv \Batrix{R}}{\phi}.
\end{align*}
The standard properties of unitarily invariant norms of operators \cite[\S~3.3]{gohberg1969introduction} provide bounds for the last two terms:
\begin{align*}
    \Norm{\Batrix{C}_{\Batrix{B}} \Batrix{G}^\pinv \Batrix{R}_{\Batrix{B}} - \Batrix{C}_{\Batrix{B}} \Batrix{G}^\pinv \Batrix{R}}{\phi} &\leq \Norm{\Batrix{C}_{\Batrix{B}} \Batrix{G}^\pinv}{\ell_2(\Hilb) \to \ell_2(\Hilb)} \Norm{\Batrix{R}_{\Batrix{E}}}{\phi}, \\
    \Norm{\Batrix{C}_{\Batrix{B}} \Batrix{G}^\pinv \Batrix{R} - \Batrix{C} \Batrix{G}^\pinv \Batrix{R}}{\phi} &\leq \Norm{\Batrix{C}_{\Batrix{E}}}{\phi} \Norm{\Batrix{G}^\pinv \Batrix{R}}{\ell_2 \to \ell_2}.
\end{align*}
To bound the first term, we refer to \cref{proposition:cgr_exact} again and rely on the equalities $\Batrix{C}_{\Batrix{B}} = \Batrix{C}_{\Batrix{B}} \Batrix{G}_{\Batrix{B}}^\pinv \Batrix{G}_{\Batrix{B}}$ and $\Batrix{R}_{\Batrix{B}} = \Batrix{G}_{\Batrix{B}} \Batrix{G}_{\Batrix{B}}^\pinv \Batrix{R}_{\Batrix{B}}$.
When combined with the defining equations of pseudoinverses, this leads to the equality
\begin{equation}
\label{eq:cgr_error_part}
    \Norm{\Batrix{C}_{\Batrix{B}} \Batrix{G}_{\Batrix{B}}^\pinv \Batrix{R}_{\Batrix{B}} - \Batrix{C}_{\Batrix{B}} \Batrix{G}^\pinv \Batrix{R}_{\Batrix{B}}}{\phi} = \Norm{\Batrix{C}_{\Batrix{B}} \Batrix{G}_{\Batrix{B}}^\pinv \Batrix{G}_{\Batrix{B}} \Batrix{G}_{\Batrix{B}}^\pinv \Batrix{R}_{\Batrix{B}} - \Batrix{C}_{\Batrix{B}} \Batrix{G}_{\Batrix{B}}^\pinv (\Batrix{G} - \Batrix{G}_{\Batrix{E}}) \Batrix{G}^\pinv (\Batrix{G} - \Batrix{G}_{\Batrix{E}}) \Batrix{G}_{\Batrix{B}}^\pinv \Batrix{R}_{\Batrix{B}}}{\phi},
\end{equation}
which we bound in turn using the triangle inequality:
\begin{align*}
    \Norm{\Batrix{C}_{\Batrix{B}} \Batrix{G}_{\Batrix{B}}^\pinv \Batrix{R}_{\Batrix{B}} - \Batrix{C}_{\Batrix{B}} \Batrix{G}^\pinv \Batrix{R}_{\Batrix{B}}}{\phi} &\leq \Norm{\Batrix{C}_{\Batrix{B}} \Batrix{G}_{\Batrix{B}}^\pinv}{\ell_2(\Hilb) \to \ell_2(\Hilb)} \Norm{\Batrix{G}_{\Batrix{B}}^\pinv \Batrix{R}_{\Batrix{B}}}{\ell_2 \to \ell_2} \\
    &\cdot \big( \Norm{\Batrix{G}_{\Batrix{E}}}{\phi} + \Norm{\Batrix{G} \Batrix{G}^\pinv \Batrix{G}_{\Batrix{E}}}{\phi} + \Norm{\Batrix{G}_{\Batrix{E}} \Batrix{G}^\pinv \Batrix{G}}{\phi} + \Norm{\Batrix{G}_{\Batrix{E}} \Batrix{G}^\pinv \Batrix{G}_{\Batrix{E}}}{\phi} \big).
\end{align*}
Since $\Batrix{G} \Batrix{G}^\pinv$ and $\Batrix{G}^\pinv \Batrix{G}$ are orthogonal projection operators in $\Hilb^{|\Index{I}|}$ and $\F^{|\Index{J}|}$, respectively, the expression in parentheses is bounded by\footnote{Owing to the normalization of symmetric gauge functions, as adopted throughout the text.} $(3 + \Norm{\Batrix{G}_{\Batrix{E}}}{\phi} \Norm{\Batrix{G}^\pinv}{\ell_2(\Hilb) \to \ell_2}) \Norm{\Batrix{G}_{\Batrix{E}}}{\phi}$.
To finish the proof, we need to bound the operator norms of $\Batrix{C}_{\Batrix{B}} \Batrix{G}_{\Batrix{B}}^\pinv$, $\Batrix{G}_{\Batrix{B}}^\pinv \Batrix{R}_{\Batrix{B}}$, $\Batrix{C}_{\Batrix{B}} \Batrix{G}^\pinv$, $\Batrix{G}^\pinv \Batrix{R}$.
The first two bounds follow from \cref{lemma:sub_pinv}.
For the third bound, consider
\begin{align*}
    \Norm{\Batrix{C}_{\Batrix{B}} \Batrix{G}^\pinv}{\ell_2(\Hilb) \to \ell_2(\Hilb)} &\leq \Norm{\Batrix{C}_{\Batrix{B}} \Batrix{G}_{\Batrix{B}}^\pinv}{\ell_2(\Hilb) \to \ell_2(\Hilb)} \Norm{\Batrix{G}_{\Batrix{B}} \Batrix{G}^\pinv}{\ell_2(\Hilb) \to \ell_2(\Hilb)} \\
    &= \Norm{\Batrix{R}_{\Batrix{U}}^\pinv}{\ell_2(\Hilb) \to \ell_2} \Norm{\Batrix{G}_{\Batrix{B}} \Batrix{G}^\pinv}{\ell_2(\Hilb) \to \ell_2(\Hilb)} \\
    &\leq \Norm{\Batrix{R}_{\Batrix{U}}^\pinv}{\ell_2(\Hilb) \to \ell_2} \big(1 + \Norm{\Batrix{G}_{\Batrix{E}} \Batrix{G}^\pinv}{\ell_2(\Hilb) \to \ell_2(\Hilb)} \big) \\
    &\leq \Norm{\Batrix{R}_{\Batrix{U}}^\pinv}{\ell_2(\Hilb) \to \ell_2} \big(1 + \Norm{\Batrix{G}_{\Batrix{E}}}{\phi} \Norm{\Batrix{G}^\pinv}{\ell_2(\Hilb) \to \ell_2} \big),
\end{align*}
where we use $\Norm{\Batrix{G}_{\Batrix{E}}}{\ell_2 \to \ell_2(\Hilb)} \leq \Norm{\Batrix{G}_{\Batrix{E}}}{\phi}$ to simplify the overall bound.
Likewise,
\begin{equation*}
    \Norm{\Batrix{G}^\pinv \Batrix{R}}{\ell_2 \to \ell_2} \leq \Norm{\Matrix{C}_{\Matrix{V}^\ast}^\pinv}{\ell_2 \to \ell_2} \big(1 + \Norm{\Batrix{G}_{\Batrix{E}}}{\phi} \Norm{\Batrix{G}^\pinv}{\ell_2(\Hilb) \to \ell_2} \big) + \Norm{\Batrix{G}^\pinv}{\ell_2(\Hilb) \to \ell_2} \Norm{\Batrix{R}_{\Batrix{E}}}{\phi}.
\end{equation*}
Combining all the intermediate bounds results in the final bound.
\end{proof}

\cref{theorem:cgr_error} quantifies how well the approximate one-sided interpolation recovers an underlying Bochner matrix $\Batrix{B}$ with low column rank; an error bound for $\Batrix{A}$ itself follows by triangle inequality.
Notably, the low-rank assumption is not imposed explicitly; rather, we replace it with a rank-equality condition, essentially demanding that the selected index sets be sufficiently rich for the particular $\Batrix{B}$.

The decomposition produced by \cref{theorem:cgr_error} possesses an inner dimension $|\Index{J}|$ and satisfies a rank bound $\Rank{c}{\Batrix{C} \Batrix{G}^\pinv \Batrix{R}} \leq \Rank{c}{\Batrix{G}} \leq |\Index{J}|$.
For a fixed selection of indices, we can further reduce the column rank by truncating $\Batrix{G}$ to lower column rank prior to taking its pseudoinverse.
Recall that $\truncate{\Batrix{G}}{k}$ denotes the SVD of $\Batrix{G}$ truncated to $k$ nonzero singular values.
Similarly, let us denote by $\truncate{\Batrix{G}}{\tau}$ the result of discarding all singular values below $\tau > 0$; we distinguish between the two based on whether the subscript is an integer or a real number.
With $\truncate{\Batrix{G}}{\tau}$ used in place of $\Batrix{G}$, the second-order error term can be controlled better, overcoming the stability issues due to the possible ill-conditioning of $\Batrix{G}$---which is common when the indices are generously oversampled.
However, the interpolation properties break down once the pseudoinverse is regularized.

\begin{corollary}
\label{corollary:cgr_error_tau}
Let $\tau > 0$.
In the setting of \cref{theorem:cgr_error},
\begin{align*}
    \Norm{\Batrix{B} - \Batrix{C} \truncate{\Batrix{G}}{\tau}^\pinv \Batrix{R}}{\phi} &\leq \eta_r \Norm{\Batrix{R}_{\Batrix{E}}}{\phi} + \eta_c \Norm{\Batrix{C}_{\Batrix{E}}}{\phi} + \eta_r \eta_c (2 \Norm{\Batrix{G}_{\Batrix{E}}}{\phi} + \Norm{\Batrix{G}_{\Batrix{B}} - \truncate{\Batrix{G}}{\tau}}{\phi}) \\
    &+ \frac{1}{\tau} \Big( \Norm{\Batrix{R}_{\Batrix{E}}}{\phi} \Norm{\Batrix{C}_{\Batrix{E}}}{\phi} + \big(\eta_r \Norm{\Batrix{R}_{\Batrix{E}}}{\phi} + \eta_c \Norm{\Batrix{C}_{\Batrix{E}}}{\phi} + \eta_r \eta_c \Norm{\Batrix{G}_{\Batrix{E}}}{\phi}\big) \Norm{\Batrix{G}_{\Batrix{E}}}{\phi} \Big).
\end{align*}
\end{corollary}
\begin{proof}
The proof of \cref{theorem:cgr_error} extends almost verbatim.
The only step that requires special attention is \Cref{eq:cgr_error_part}.
Expanding its right-hand side, we note that $\Batrix{G} \truncate{\Batrix{G}}{\tau}^\pinv \Batrix{G} = \truncate{\Batrix{G}}{\tau}^\pinv$ and, subsequently, rely on $\Batrix{G} \truncate{\Batrix{G}}{\tau}^\pinv$ and $\truncate{\Batrix{G}}{\tau}^\pinv \Batrix{G}$ also being orthogonal projection operators due to the SVD truncation.
\end{proof}

A similar corollary holds if we truncate the rank rather than apply thresholding to the singular values.
Below, we bound the norm of the pseudoinverse of the truncated $\Batrix{G}$ in terms of $\Batrix{G}_{\Batrix{B}}$ and $\Batrix{G}_{\Batrix{E}}$, separating their contributions.

\begin{corollary}
\label{corollary:cgr_error_k}
In the setting of \cref{theorem:cgr_error}, let $1 \leq k \leq \Rank{c}{\Batrix{G}_{\Batrix{B}}}$ be such that $\gamma = \sigma_{k}(\Batrix{G}_{\Batrix{B}}) - \Norm{\Batrix{G}_{\Batrix{E}}}{\ell_2 \to \ell_2(\Hilb)} > 0$.
Then
\begin{align*}
    \Norm{\Batrix{B} - \Batrix{C} \truncate{\Batrix{G}}{k}^\pinv \Batrix{R}}{\phi} &\leq \eta_r \Norm{\Batrix{R}_{\Batrix{E}}}{\phi} + \eta_c \Norm{\Batrix{C}_{\Batrix{E}}}{\phi} + \eta_r \eta_c (2 \Norm{\Batrix{G}_{\Batrix{E}}}{\phi} + \Norm{\Batrix{G}_{\Batrix{B}} - \truncate{\Batrix{G}}{k}}{\phi}) \\
    &+ \frac{1}{\gamma} \Big( \Norm{\Batrix{R}_{\Batrix{E}}}{\phi} \Norm{\Batrix{C}_{\Batrix{E}}}{\phi} + \big(\eta_r \Norm{\Batrix{R}_{\Batrix{E}}}{\phi} + \eta_c \Norm{\Batrix{C}_{\Batrix{E}}}{\phi} + \eta_r \eta_c \Norm{\Batrix{G}_{\Batrix{E}}}{\phi}\big) \Norm{\Batrix{G}_{\Batrix{E}}}{\phi} \Big).
\end{align*}
\end{corollary}
\begin{proof}
By Weyl's inequalities, which apply to Bochner matrices, $\sigma_k(\Batrix{G}) \geq \gamma > 0$, and thus the rank truncation $\truncate{\Batrix{G}}{k}$ is well-defined.
The proof remains the same, except for the bound $\Norm{\truncate{\Batrix{G}}{k}^\pinv}{\ell_2(\Hilb) \to \ell_2} = \sigma_k(\Batrix{G})^{-1} \leq \gamma^{-1}$.
\end{proof}

\subsection{Well-conditioned submatrices}
\label{subsec:rank_excess}
The approximation error bound in \cref{theorem:cgr_error} is governed by the operator norms of the pseudoinverses of submatrices of singular factors.
The existence of well-conditioned submatrices of classical $n \times \kappa$ matrices with orthonormal columns is a long-standing question.
It was conjectured in \cite{goreinov1997theory} that there always exists a $\kappa \times \kappa$ submatrix for which the operator norm of its inverse is bounded by $\sqrt{n}$; see \cite{nesterenko2025subspaces, sengupta2026submatrices} for recent progress in this direction.

\emph{Maximum-volume} submatrices do not resolve this conjecture, but are nonetheless well-conditioned \cite{goreinov1997theory, de2007subset, mikhalev2018rectangular, osinsky2018pseudo}.
The volume is defined as $\nu(\Matrix{A}) = \sqrt{\det(\Matrix{A}^\ast \Matrix{A})}$, and if a $|\Index{J}| \times \kappa$ ($|\Index{J}| \geq \kappa$) submatrix has the largest volume among all $|\Index{J}| \times \kappa$ submatrices then the operator norm of its pseudoinverse is upper bounded by
\begin{equation}
\label{eq:maxvol_classical}
    \sqrt{1 + \frac{(n - |\Index{J}|)\kappa}{|\Index{J}| - \kappa + 1}}.   
\end{equation}
We generalize this to Bochner matrices and define their volume as $\nu(\Batrix{A}) = \sqrt{\det(\Batrix{A}^\ast \Batrix{A})}$, stressing that $\nu(\Batrix{A}^\trans) \neq \nu(\Batrix{A})$ in general. 
In the following statements, we also denote by $\Batrix{A}_{-i}$ the Bochner matrix with its $i$th row $\Batrix{a}_i^\trans$ removed.

\begin{lemma}
\label{lemma:rank_excess}
Let $\Batrix{A} \in \Hilb^{N \times \kappa}$ have $\Rank{c}{\Batrix{A}} = \kappa$.
Let $\Matrix{\Gamma} = \Batrix{A}^\ast \Batrix{A}$ and $\Matrix{\Gamma}_i = (\Batrix{a}_{i}^\trans)^\ast \Batrix{a}_{i}^\trans$ for every $1 \leq i \leq N$.
Then
\begin{equation*}
    \frac{\nu^2(\Batrix{A}_{-i})}{\nu^2(\Batrix{A})} = \det(\Matrix{I}_{\kappa} - \Matrix{S}_i), \quad \Matrix{S}_i = \Matrix{\Gamma}^{-\tfrac{1}{2}} \Matrix{\Gamma}_i \Matrix{\Gamma}^{-\tfrac{1}{2}}.
\end{equation*}
Let $\Delta_{i}(\Batrix{A}) = \det(\Matrix{I}_{\kappa} - \Matrix{S}_i) - (1 - \trace(\Matrix{S}_i))$ and $\Delta(\Batrix{A}) = \sum_{i = 1}^{N} \Delta_{i}(\Batrix{A})$.
Then
\begin{equation*}
    \sum_{i = 1}^{N} \frac{\nu^2(\Batrix{A}_{-i})}{\nu^2(\Batrix{A})} = N - \kappa + \Delta(\Batrix{A}).
\end{equation*}
\end{lemma}
\begin{proof}
By assumption, $\Matrix{\Gamma} \in \F^{\kappa \times \kappa}$ has full rank and $\nu(\Batrix{A}) \neq 0$; consequently, all $\Matrix{S}_i$ are well-defined positive semidefinite matrices.
Furthermore, $\Matrix{\Gamma} = \sum_{i = 1}^{N} \Matrix{\Gamma}_i$ and $\sum_{i = 1}^{N} \Matrix{S}_i = \Matrix{I}_{\kappa}$.
We expand the definition of $\nu^2(\Batrix{A}_{-i})$ to prove the first equality:
\begin{equation*}
    \nu^2(\Batrix{A}_{-i}) = \det(\Batrix{A}_{-i}^\ast \Batrix{A}_{-i}) = \det(\Matrix{\Gamma} - \Matrix{\Gamma}_i) = \nu^2(\Batrix{A}) \det(\Matrix{I}_{\kappa} - \Matrix{S}_i).
\end{equation*}
To obtain the second equality, observe that
\begin{equation*}
    \sum_{i = 1}^{N} (1 - \trace(\Matrix{S}_i)) = N - \trace\left(\sum_{i = 1}^N \Matrix{S}_i\right) = N - \trace(\Matrix{I}_{\kappa}) = N - \kappa.
\end{equation*}
\end{proof}

We call $\Delta_i(\Batrix{A})$ the \emph{intra-row rank excess} of $\Batrix{A}$ and $\Delta(\Batrix{A})$ its \emph{inter-row rank excess}.
Since the matrices $\Matrix{S}_i$ are positive semidefinite and satisfy $\sum_{i = 1}^{N} \Matrix{S}_i = \Matrix{I}_{\kappa}$, their eigenvalues lie in $[0, 1]$.
Thus, by the Weierstrass product inequality, we have $\Delta_i(\Batrix{A}) \geq 0$; moreover, $\Delta_i(\Batrix{A})$ vanishes if and only if $\rank{\Matrix{\Gamma}_i} \leq 1$.
Based on the SVD of $\Batrix{a}_i^\trans$, we can establish a simple yet fundamental relationship:
\begin{equation*}
    \rank{\Matrix{\Gamma}_i} = \Rank{c}{\Batrix{a}_i^\trans}.
\end{equation*}
Therefore, the inter-row rank excess is always zero for classical matrices.

\begin{lemma}
\label{lemma:vol_ratio_pinv_bound}
In the setting of \cref{lemma:rank_excess}, let $i_\ast$ be such that $\Rank{c}{\Batrix{A}_{-i_\ast}} = \kappa$, and denote $\Matrix{\Gamma}_{-i_\ast} = \Batrix{A}_{-i_\ast}^\ast \Batrix{A}_{-i_\ast}$.
Then
\begin{equation*}
    \frac{\nu^2(\Batrix{A})}{\nu^2(\Batrix{A}_{-i_\ast})} = \det\Big(\Matrix{I}_{\kappa} + \Matrix{\Gamma}_{-i_\ast}^{-\tfrac{1}{2}} \Matrix{\Gamma}_{i_\ast} \Matrix{\Gamma}_{-i_\ast}^{-\tfrac{1}{2}}\Big) \geq 1 + \Norm{\Batrix{a}_{i_\ast}^\trans \Batrix{A}_{-i_\ast}^{\pinv}}{\ell_2(\Hilb) \to \Hilb}^2.
\end{equation*}
\end{lemma}
\begin{proof}
The equality follows from the definition of volume.
To prove the inequality, consider vectors $\Matrix{x} \in \F^{\kappa}$ of unit $\ell_2$ norm and let $\Batrix{A}_{-i_\ast} = \Batrix{P} \Matrix{\Lambda} \Matrix{Q}^\ast$ be an SVD.
Then
\begin{align*}
    \Norm[\Big]{\Matrix{\Gamma}_{-i_\ast}^{-\tfrac{1}{2}} \Matrix{\Gamma}_{i_\ast} \Matrix{\Gamma}_{-i_\ast}^{-\tfrac{1}{2}}}{\ell_2 \to \ell_2} &= \max_{\Matrix{x}} \Big|
    \Matrix{x}^\ast \Matrix{\Gamma}_{-i_\ast}^{-\tfrac{1}{2}} \Matrix{\Gamma}_{i_\ast} \Matrix{\Gamma}_{-i_\ast}^{-\tfrac{1}{2}} \Matrix{x} \Big| = \max_{\Matrix{x}} \Norm[\Big]{\Batrix{a}_{i_\ast}^\trans \Matrix{\Gamma}_{-i_\ast}^{-\tfrac{1}{2}} \Matrix{x}}{\Hilb}^2 = \Norm[\Big]{\Batrix{a}_{i_\ast}^\trans \Matrix{\Gamma}_{-i_\ast}^{-\tfrac{1}{2}}}{\ell_2 \to \Hilb}^2 \\
    &= \Norm[\Big]{\Batrix{a}_{i_\ast}^\trans \Matrix{Q} \Matrix{\Lambda}^{-1}}{\ell_2 \to \Hilb}^2 = \Norm[\Big]{\Batrix{a}_{i_\ast}^\trans \Matrix{Q} \Matrix{\Lambda}^{-1} \Batrix{P}^\ast}{\ell_2(\Hilb) \to \Hilb}^2 = \Norm{\Batrix{a}_{i_\ast}^\trans \Batrix{A}_{-i_\ast}^{\pinv}}{\ell_2(\Hilb) \to \Hilb}^2,
\end{align*}
where we use the unitary invariance of operator norms.
The inequality follows since the remaining eigenvalues are nonnegative.
\end{proof}

Combining \cref{lemma:rank_excess,lemma:vol_ratio_pinv_bound}, we can bound the volume of $\Batrix{A}$ as
\begin{equation*}
    \Big( 1 + \Norm{\Batrix{a}_{i_\ast}^\trans \Batrix{A}_{-i_\ast}^{\pinv}}{\ell_2(\Hilb) \to \Hilb}^2 \Big) \nu^2(\Batrix{A}_{-i_\ast}) \leq \nu^2(\Batrix{A}) = \frac{N}{N - \kappa + \Delta(\Batrix{A})} \bigg(\frac{1}{N} \sum_{i = 1}^{N} \nu^2(\Batrix{A}_{-i})\bigg).
\end{equation*}
Suppose that the squared volume of the submatrix $\Batrix{A}_{-i_\ast}$ is at least the average squared volume across all $(N-1) \times \kappa$ submatrices of $\Batrix{A}$.
Then this volume bound yields
\begin{equation}
\label{eq:maxvol_lsq_bound}
    \Norm{\Batrix{a}_{i_\ast}^\trans \Batrix{A}_{-i_\ast}^{\pinv}}{\ell_2(\Hilb) \to \Hilb}^2 \leq \frac{\kappa - \Delta(\Batrix{A})}{N - \kappa + \Delta(\Batrix{A})},
\end{equation}
which extends the classical maximum-volume least-squares bound to Bochner matrices and relies on a relaxed above-average volume condition.
We are finally ready to derive an operator-norm bound for the pseudoinverse of a maximum-volume submatrix.

\begin{theorem}
\label{theorem:maxvol_pinv_bound}
Let $\Batrix{U} \in \Hilb^{m \times \kappa}$ have orthonormal columns.
Let indices $\Index{I} \subseteq [m]$ be such that $\Batrix{R}_{\Batrix{U}}$ has the largest volume among all $|\Index{I}| \times \kappa$ submatrices of $\Batrix{U}$.
If $\nu(\Batrix{R}_{\Batrix{U}}) > 0$,
\begin{equation*}
    \Norm{\Batrix{R}_{\Batrix{U}}^\pinv}{\ell_2(\Hilb) \to \ell_2} \leq \sqrt{1 + \sum_{i \not\in \Index{I}} \frac{\kappa - \Delta(\tilde{\Batrix{R}}_{i})}{|\Index{I}| + 1 - \kappa + \Delta(\tilde{\Batrix{R}}_{i})}}, \quad \tilde{\Batrix{R}}_{i} = \begin{bmatrix}
        \Batrix{R}_{\Batrix{U}} \\
        \Batrix{u}_{i}^\trans
    \end{bmatrix}.
\end{equation*}
\end{theorem}
\begin{proof}
Using the standard properties of the operator norm, we have
\begin{align*}
    \Norm{\Batrix{R}_{\Batrix{U}}^\pinv}{\ell_2(\Hilb) \to \ell_2}^2 &= \Norm{\Batrix{U} \Batrix{R}_{\Batrix{U}}^\pinv}{\ell_2(\Hilb) \to \ell_2(\Hilb)}^2 = \Norm{(\Batrix{U} \Batrix{R}_{\Batrix{U}}^\pinv)^\ast (\Batrix{U} \Batrix{R}_{\Batrix{U}}^\pinv)}{\ell_2(\Hilb) \to \ell_2(\Hilb)} \\
    &= \Norm[\bigg]{(\Batrix{R}_{\Batrix{U}} \Batrix{R}_{\Batrix{U}}^\pinv)^\ast (\Batrix{R}_{\Batrix{U}} \Batrix{R}_{\Batrix{U}}^\pinv) + \sum_{i \not\in \Index{I}} (\Batrix{u}_{i}^\trans \Batrix{R}_{\Batrix{U}}^\pinv)^\ast (\Batrix{u}_{i}^\trans \Batrix{R}_{\Batrix{U}}^\pinv)}{\ell_2(\Hilb) \to \ell_2(\Hilb)} \\
    &= \Norm[\bigg]{\Batrix{R}_{\Batrix{U}} \Batrix{R}_{\Batrix{U}}^\pinv + \sum_{i \not\in \Index{I}} (\Batrix{u}_{i}^\trans \Batrix{R}_{\Batrix{U}}^\pinv)^\ast (\Batrix{u}_{i}^\trans \Batrix{R}_{\Batrix{U}}^\pinv)}{\ell_2(\Hilb) \to \ell_2(\Hilb)}\\
    &\leq 1 + \sum_{i \not\in \Index{I}} \Norm{(\Batrix{u}_{i}^\trans \Batrix{R}_{\Batrix{U}}^\pinv)^\ast (\Batrix{u}_{i}^\trans \Batrix{R}_{\Batrix{U}}^\pinv)}{\ell_2(\Hilb) \to \ell_2(\Hilb)} \\
    &= 1 + \sum_{i \not\in \Index{I}} \Norm{\Batrix{u}_{i}^\trans \Batrix{R}_{\Batrix{U}}^\pinv}{\ell_2(\Hilb) \to \Hilb}^2,
\end{align*}
since $\Batrix{R}_{\Batrix{U}} \Batrix{R}_{\Batrix{U}}^\pinv$ is an orthogonal projector.
We then bound each term in the sum using \cref{eq:maxvol_lsq_bound}; this is allowed because the positive volume of $\Batrix{R}_{\Batrix{U}}$ implies $\Rank{c}{\Batrix{R}_{\Batrix{U}}} = \kappa$.
\end{proof}

For classical matrices, \cref{theorem:maxvol_pinv_bound} recovers the standard well-conditioning bound for maximum-volume submatrices.
Meanwhile, the presence of inter-row rank excesses in the bound shows that maximum-volume submatrices are \emph{even better conditioned} in the Bochner setting; indeed, the right-hand side in \cref{eq:maxvol_lsq_bound} is a monotonically decreasing function of the inter-row rank excess.

To illustrate the conditioning gains, consider $\Batrix{Q} \in \Hilb^{m \times \kappa}$ whose entries are mutually orthonormal in $\Hilb$ and set $\Batrix{U} = \Batrix{Q} / \sqrt{m}$.
Due to symmetry, all $|\Index{I}| \times \kappa$ submatrices have the same volume, and it holds for every row-augmentation in \cref{theorem:maxvol_pinv_bound} that
\begin{equation*}
    \Delta(\tilde{\Batrix{R}}_{i}) = (|\Index{I}| + 1) \left( 1 - \frac{1}{|\Index{I}| + 1} \right)^{\kappa} - (|\Index{I}| + 1 - \kappa),
\end{equation*}
which we can substitute into the bound, simplifying to
\begin{equation*}
     \Norm{\Batrix{R}_{\Batrix{U}}^\pinv}{\ell_2(\Hilb) \to \ell_2} \leq \sqrt{1 + (m - |\Index{I}|) \Big[ \Big(1 + \frac{1}{|\Index{I}|}\Big)^\kappa - 1 \Big]} < \sqrt{1 + (m - |\Index{I}|) \big( e^{\kappa / |\Index{I}|} - 1 \big)}.
\end{equation*}
It is instructive to compare this bound with the bound for classical maximum-volume submatrices \cref{eq:maxvol_classical}; see \cref{tab:maxvol_comparison}.
When a square submatrix is selected, classical matrices enjoy the bound $\sqrt{1 + m \kappa}$, whereas the bound for $\Batrix{U}$ is $\sqrt{1 + (e-1)m}$ and does not grow with $\kappa$.
A standard way to regularize classical submatrices is to sample twice as many rows as there are columns, improving the bound to $\sqrt{1 + m}$.
For the specific $\Batrix{U}$, such oversampling leads to an even smaller bound $\sqrt{1 + (e^{1/2}-1)m}$.
Unlike classical matrices, the rows of Bochner matrices can be \emph{undersampled}, while still delivering well-conditioned submatrices.
For example, the selection of $\kappa / \log(\kappa)$ rows results in a bound that is similar to the case of classical square submatrices.

\begin{table}[htb!]
\centering
\caption{Comparison of operator-norm bounds for the pseudoinverse of maximum-volume submatrices of classical and Bochner matrices with orthonormal columns.}
\label{tab:maxvol_comparison}
\begin{tabular*}{\textwidth}{@{\extracolsep\fill}lcc@{}}
    \toprule
    $|\Index{I}|$ & Classical matrices & Bochner matrix $\Batrix{U}$ \\
    \midrule
    $\kappa / \log(\kappa)$  & --- & $\sqrt{1 + (m - \kappa / \log(\kappa)) (\kappa - 1)}$ \\
    $\kappa$  & $\sqrt{1 + (m - \kappa) \kappa}$ & $\sqrt{1 + (m - \kappa) (e - 1)}$ \\
    $c\kappa,~c \geq 1$  & $\sqrt{1 + (m - c\kappa) \frac{\kappa}{(c-1)\kappa + 1}}$ & $\sqrt{1 + (m - c\kappa) (e^{1/c} - 1)}$ \\
    \bottomrule
\end{tabular*}
\end{table}

Let us note that \cref{theorem:maxvol_pinv_bound} holds for any $\Batrix{R}_{\Batrix{U}}$ that has the largest volume among all submatrices that disagree with it in exactly one row.
Relaxed further, it is sufficient for $\Batrix{R}_{\Batrix{U}}$ to have above-average volume within each individual augmentation $\tilde{\Batrix{R}}_{i}$.

\subsection{Entrywise error bound for one-sided interpolation}
Another important class of approximation guarantees for classical cross approximation are entrywise error bounds for when the submatrix at the intersection of the selected columns and rows itself exhibits maximum-volume properties.
Following \cite{osinsky2018pseudo}, we define the \emph{projective volume} of a Bochner matrix as $\nu_{\kappa}(\Batrix{A}) = \prod_{s = 1}^{\kappa} \sigma_s(\Batrix{A})$ and truncate the rank of the pseudoinverse for regularization.
We also introduce a \emph{projective inter-row rank excess} $\Delta_{\kappa}(\Batrix{A}) = \sup_{\Matrix{Q}_{\kappa}} \Delta(\Batrix{A} \Matrix{Q}_{\kappa})$, where $\Matrix{Q}_\kappa$ are (non-unique) truncated right singular factors of $\Batrix{A}$.
This quantity is well-defined when $\Rank{c}{\Batrix{A}} \geq \kappa$.

\begin{theorem}
\label{theorem:cgr_error_maxvol}
Let $\Batrix{A} \in \Hilb^{m \times n}$ and $\kappa \in \N$.
Let $\Batrix{G}$ have the maximum $\kappa$-projective volume among all $|\Index{I}| \times |\Index{J}|$ submatrices of $\Batrix{A}$.
If $\nu_{\kappa}(\Batrix{G}) > 0$,
\begin{equation*}
    \Norm{\Batrix{A} - \Batrix{C} \truncate{\Batrix{G}}{\kappa}^\pinv \Batrix{R}}{\ell_{\infty}(\Hilb)} \leq \sqrt{1 + \frac{\kappa - \Delta_\kappa'}{|\Index{I}| + 1 - \kappa + \Delta_\kappa'}} \sqrt{1 + \frac{\kappa}{|\Index{J}| + 1 - \kappa}} \sigma_{\kappa + 1}(\Batrix{A})
\end{equation*}
with $\Delta_\kappa' = \min_{i \not\in \Index{I}} \min_{1 \leq j \leq n} \Delta_{\kappa}\big( \Batrix{A}(\Index{I} \cup \{i\}, \Index{J} \cup \{j\}) \big)$.
\end{theorem}
\begin{proof}
For fixed indices $1 \leq i \leq m$ and $1 \leq j \leq n$, denote
\begin{equation*}
    \Bentry{a} = \Batrix{A}(i, j), \quad \Batrix{c}^\trans = \Batrix{C}(i, :), \quad  \Batrix{r} = \Batrix{R}(:, j), \quad \Bentry{e} = \Batrix{a} - \Batrix{c}^\trans \truncate{\Batrix{G}}{\kappa}^\pinv \Batrix{r}. 
\end{equation*}
Our goal is to bound $\Norm{\Bentry{e}}{\Hilb}$ depending on whether $i$ and $j$ fall into the index sets $\Index{I}$ and $\Index{J}$, respectively.
We shall assume that $\Bentry{e} \neq 0$, for otherwise the bounds are trivial.

First, let $i \in \Index{I}$ and $j \in \Index{J}$.
Then
\begin{align*}
    \Norm{\Bentry{e}}{\Hilb} &\leq \Norm{\Batrix{G} - \Batrix{G} \truncate{\Batrix{G}}{\kappa}^\pinv \Batrix{G}}{\ell_{\infty}(\Hilb)} = \Norm{\Batrix{G} - \truncate{\Batrix{G}}{\kappa}}{\ell_{\infty}(\Hilb)} \\
    &\leq \Norm{\Batrix{G} - \truncate{\Batrix{G}}{\kappa}}{\ell_2 \to \ell_2(\Hilb)} = \sigma_{\kappa + 1}(\Batrix{G}) \leq \sigma_{\kappa + 1}(\Batrix{A}),
\end{align*}
where we use standard properties of the pseudoinverse and singular values.

Second, let $i \not\in \Index{I}$ and $j \not\in \Index{J}$.
Construct an augmented submatrix
\begin{equation*}
    \Batrix{M} = \begin{bmatrix}
        \Batrix{G} & \Batrix{r} \\
        \Batrix{c}^\trans & \Batrix{a}
    \end{bmatrix} \in \Hilb^{(|\Index{I}| + 1) \times (|\Index{J}| + 1)}
\end{equation*}
and let $\truncate{\Batrix{G}}{\kappa} = \Batrix{U}_{\kappa} \Matrix{\Sigma}_{\kappa} \Matrix{V}_{\kappa}^\ast$ be a truncated SVD.
Consider an auxiliary matrix
\begin{equation*}
    \Matrix{T} = \begin{bmatrix}
        \Matrix{V}_{\kappa} & -\truncate{\Batrix{G}}{\kappa}^\pinv \Batrix{r} \\
        0 & 1
    \end{bmatrix} \in \F^{(|\Index{J}| + 1) \times (\kappa + 1)}
\end{equation*}
and multiply $\Batrix{M}$ with it on the right:
\begin{equation*}
    \Batrix{K} = \Batrix{M} \Matrix{T} = \begin{bmatrix}
        \Batrix{U}_{\kappa} \Matrix{\Sigma}_{\kappa} & \Batrix{r}^\perp \\
        \Batrix{c}^\trans \Matrix{V}_{\kappa} & \Bentry{e}
    \end{bmatrix} \in \Hilb^{(|\Index{I}| + 1) \times (\kappa + 1)}, \quad \Batrix{r}^\perp = \Batrix{r} - \Batrix{U}_{\kappa} \Batrix{U}_{\kappa}^\ast \Batrix{r} \in \Hilb^{|\Index{I}|}.
\end{equation*}
By the definition of volume, we get
\begin{equation*}
    \nu^2(\Batrix{K}) = \det(\Batrix{K}^\ast \Batrix{K}) = \det \left( \begin{bmatrix}
        \Matrix{\Sigma}_{\kappa}^2 + (\Batrix{c}^\trans \Matrix{V}_{\kappa})^\ast \Batrix{c}^\trans \Matrix{V}_{\kappa} & (\Batrix{c}^\trans \Matrix{V}_{\kappa})^\ast \Bentry{e} \\
        \Bentry{e}^\ast \Batrix{c}^\trans \Matrix{V}_{\kappa} & \Norm{\Batrix{r}^\perp}{\ell_2(\Hilb)}^2 + \Norm{\Bentry{e}}{\Hilb}^2
    \end{bmatrix} \right).
\end{equation*}
Since $\Bentry{e} \neq 0$ by our assumption, we can invoke Schur's complement:
\begin{equation*}
    \nu^2(\Batrix{K}) = (\Norm{\Batrix{r}^\perp}{\ell_2(\Hilb)}^2 + \Norm{\Bentry{e}}{\Hilb}^2) \det\bigg( \Matrix{\Sigma}_{\kappa}^2 + (\Batrix{c}^\trans \Matrix{V}_{\kappa})^\ast \Batrix{c}^\trans \Matrix{V}_{\kappa} - \frac{(\Batrix{c}^\trans \Matrix{V}_{\kappa})^\ast \Bentry{e} \Bentry{e}^\ast \Batrix{c}^\trans \Matrix{V}_{\kappa}}{\Norm{\Batrix{r}^\perp}{\ell_2(\Hilb)}^2 + \Norm{\Bentry{e}}{\Hilb}^2} \bigg).
\end{equation*}
The matrix
\begin{equation*}
    (\Batrix{c}^\trans \Matrix{V}_{\kappa})^\ast \Batrix{c}^\trans \Matrix{V}_{\kappa} - \frac{(\Batrix{c}^\trans \Matrix{V}_{\kappa})^\ast \Bentry{e} \Bentry{e}^\ast \Batrix{c}^\trans \Matrix{V}_{\kappa}}{\Norm{\Bentry{e}}{\Hilb}^2} = (\Batrix{c}^\trans \Matrix{V}_{\kappa} - \Bentry{e} \Bentry{e}^\pinv \Batrix{c}^\trans \Matrix{V}_{\kappa})^\ast (\Batrix{c}^\trans \Matrix{V}_{\kappa} - \Bentry{e} \Bentry{e}^\pinv \Batrix{c}^\trans \Matrix{V}_{\kappa})
\end{equation*}
is a Gram matrix and hence PSD; therefore,
\begin{equation*}
    (\Batrix{c}^\trans \Matrix{V}_{\kappa})^\ast \Batrix{c}^\trans \Matrix{V}_{\kappa} - \frac{(\Batrix{c}^\trans \Matrix{V}_{\kappa})^\ast \Bentry{e} \Bentry{e}^\ast \Batrix{c}^\trans \Matrix{V}_{\kappa}}{\Norm{\Batrix{r}^\perp}{\ell_2(\Hilb)}^2 + \Norm{\Bentry{e}}{\Hilb}^2} \succcurlyeq (\Batrix{c}^\trans \Matrix{V}_{\kappa})^\ast \Batrix{c}^\trans \Matrix{V}_{\kappa} - \frac{(\Batrix{c}^\trans \Matrix{V}_{\kappa})^\ast \Bentry{e} \Bentry{e}^\ast \Batrix{c}^\trans \Matrix{V}_{\kappa}}{\Norm{\Bentry{e}}{\Hilb}^2} \succcurlyeq 0,
\end{equation*}
and the monotonicity of determinant for PSD matrices implies
\begin{equation*}
    \nu^2(\Batrix{K}) \geq (\Norm{\Batrix{r}^\perp}{\ell_2(\Hilb)}^2 + \Norm{\Bentry{e}}{\Hilb}^2) \det(\Matrix{\Sigma}_{\kappa}^2) \geq \Norm{\Bentry{e}}{\Hilb}^2 \det(\Matrix{\Sigma}_{\kappa}^2).
\end{equation*}
Let us relate the volumes of $\Batrix{K}$ and $\Batrix{M}$.
Using Schur's complement, we can show that the auxiliary matrix $\Matrix{T}$ satisfies
\begin{equation*}
    \nu_{\kappa + 1}^2(\Matrix{T}) = \nu^2(\Matrix{T}) =  \det(\Matrix{T}^\ast \Matrix{T}) = \det\left(\begin{bmatrix}
        \Matrix{I}_{\kappa} & -\Matrix{\Sigma}_\kappa^{-1} \Batrix{U}_\kappa^\ast \Batrix{r} \\
        (-\Matrix{\Sigma}_\kappa^{-1} \Batrix{U}_\kappa^\ast \Batrix{r})^\ast & 1 + \Norm{\Matrix{\Sigma}_\kappa^{-1} \Batrix{U}_\kappa^\ast \Batrix{r}}{\ell_2}^2
    \end{bmatrix}\right) = 1.
\end{equation*}
Let $\Batrix{M} = \Batrix{P} \Matrix{\Lambda} \Matrix{Q}^\ast$ be an SVD, then by the log-majorization of singular values of a matrix product \cite[Thm.~3.3.14]{horn1994topics} we obtain
\begin{equation*}
    \nu(\Batrix{K}) = \nu(\Batrix{P} \Batrix{P}^\ast \Batrix{M} \Matrix{T}) = \nu(\Batrix{P}^\ast \Batrix{M} \Matrix{T}) = \nu_{\kappa + 1}(\Batrix{P}^\ast \Batrix{M} \Matrix{T}) \leq \nu_{\kappa + 1}(\Batrix{P}^\ast \Batrix{M}) \nu_{\kappa + 1}(\Matrix{T}) = \nu_{\kappa + 1}(\Batrix{M}),
\end{equation*}
and therefore
\begin{equation*}
    \Norm{\Bentry{e}}{\Hilb} \leq \frac{\nu(\Batrix{K})}{\nu_{\kappa}(\Batrix{G})} \leq \frac{\nu_{\kappa + 1}(\Batrix{M})}{\nu_{\kappa}(\Batrix{G})} = \frac{\nu_{\kappa}(\Batrix{M})}{\nu_{\kappa}(\Batrix{G})} \sigma_{\kappa + 1}(\Batrix{M}) \leq \frac{\nu_{\kappa}(\Batrix{M})}{\nu_{\kappa}(\Batrix{G})} \sigma_{\kappa + 1}(\Batrix{A}).
\end{equation*}
To bound the volume ratio, we shall expand $\nu_{\kappa}^2(\Batrix{M})$ into a sum of squared volumes of submatrices with removed columns and rows based on the fundamental identities
\begin{equation*}
    \nu_{\kappa}(\Batrix{M}) = \nu(\Batrix{P}_{\kappa}^\ast \Batrix{M}) = \nu(\Batrix{M} \Matrix{Q}_{\kappa}) = \max_{\Batrix{\tilde{P}}^\ast \Batrix{\tilde{P}} = \Matrix{I}_{\kappa}} \nu(\Batrix{\tilde{P}}^\ast \Batrix{M}) = \max_{\Matrix{\tilde{Q}}^\ast \Matrix{\tilde{Q}} = \Matrix{I}_{\kappa}} \nu(\Batrix{M} \Matrix{\tilde{Q}}).
\end{equation*}
Denote by $\Batrix{M}_{-s}$ the submatrix with $s$th row removed, by $\Batrix{M}_{:, -t}$ the submatrix with $t$th column removed, and by $\Batrix{M}_{-s,-t}$ the submatrix with both $s$th row and $t$th column removed.
Applying \cref{lemma:rank_excess} to $\Batrix{M} \Matrix{Q}_\kappa$, which is allowed since $\nu_{\kappa}(\Batrix{G}) > 0$, we obtain
\begin{equation*}
    \sum_{s = 1}^{|\Index{I}| + 1} \nu^2(\Batrix{M}_{-s} \Matrix{Q}_\kappa) = \Big(|\Index{I}| + 1 - \kappa + \Delta(\Batrix{M} \Matrix{Q}_\kappa) \Big) \nu_{\kappa}^2(\Batrix{M}).
\end{equation*}
Next, we bound each term via $\nu(\Batrix{M}_{-s} \Matrix{Q}_\kappa) \leq \nu_{\kappa}(\Batrix{M}_{-s})$.
If $\nu_{\kappa}(\Batrix{M}_{-s}) = 0$ then we move on; otherwise, let $\Batrix{P}_{s,\kappa}$ be the truncated left singular factor of $\Batrix{M}_{-s}$.
Then we use the equality $\nu_{\kappa}(\Batrix{M}_{-s}) = \nu(\Batrix{P}_{s,\kappa}^\ast \Batrix{M}_{-s})$ and apply \cref{lemma:rank_excess} to $(\Batrix{P}_{s,\kappa}^\ast \Batrix{M}_{-s})^\trans$.
Because this is a classical matrix, its inter-row excess is zero, and hence
\begin{equation*}
    \sum_{t = 1}^{|\Index{J}| + 1} \nu^2(\Batrix{P}_{s,\kappa}^\ast \Batrix{M}_{-s,-t}) = \Big(|\Index{J}| + 1 - \kappa \Big) \nu_{\kappa}^2(\Batrix{M}_{-s}).
\end{equation*}
Combining the two sums yields
\begin{align*}
    \sum_{s = 1}^{|\Index{I}| + 1} \sum_{t = 1}^{|\Index{J}| + 1} \nu_{\kappa}^2(\Batrix{M}_{-s,-t}) &\geq \sum_{s = 1}^{|\Index{I}| + 1} \sum_{t = 1}^{|\Index{J}| + 1} \nu^2(\Batrix{P}_{s,\kappa}^\ast \Batrix{M}_{-s,-t}) = \sum_{s = 1}^{|\Index{I}| + 1} \Big(|\Index{J}| + 1 - \kappa \Big) \nu_{\kappa}^2(\Batrix{M}_{-s}) \\
    &\geq \Big(|\Index{J}| + 1 - \kappa \Big) \sum_{s = 1}^{|\Index{I}| + 1} \nu^2(\Batrix{M}_{-s} \Matrix{Q}_\kappa) \\
    &= \Big(|\Index{J}| + 1 - \kappa \Big) \Big(|\Index{I}| + 1 - \kappa + \Delta(\Batrix{M} \Matrix{Q}_\kappa) \Big) \nu_{\kappa}^2(\Batrix{M}),
\end{align*}
and it remains to recall that $\nu_{\kappa}(\Batrix{M}_{-s,-t}) \leq \nu_{\kappa}(\Batrix{G})$ and maximize over the SVDs of $\Batrix{M}$.

Third, let $i \in \Index{I}$ and $j \not\in \Index{J}$.
Then we augment as $\Batrix{M} = [\Batrix{G}~\Batrix{r}]$ and use the same auxiliary matrix $\Matrix{T}$ to form $\Batrix{K}$.
We can then readily show that 
\begin{equation*}
    \nu^2(\Batrix{K}) = \Norm{\Batrix{r}^\perp}{\ell_2(\Hilb)}^2 \det(\Matrix{\Sigma}_{\kappa}^2) \geq \Norm{\Bentry{e}}{\Hilb}^2 \det(\Matrix{\Sigma}_{\kappa}^2),
\end{equation*}
since $\Bentry{e}$ is now a component of $\Batrix{r}^\perp$, proceed to the same bound $\Norm{\Bentry{e}}{\Hilb} \leq \tfrac{\nu_{\kappa}(\Batrix{M})}{\nu_{\kappa}(\Batrix{G})} \sigma_{\kappa + 1}(\Batrix{A})$, and skip the row-removal step in volume expansions to get
\begin{equation*}
    \Norm{\Bentry{e}}{\Hilb} \leq \sqrt{1 + \frac{\kappa}{|\Index{J}| + 1 - \kappa}} \sigma_{\kappa + 1}(\Batrix{A}).
\end{equation*}

Fourth, let $i \not\in \Index{I}$ and $j \in \Index{J}$.
Then we augment as $\Batrix{M} = [\Batrix{G}^\trans~\Batrix{c}]^\trans \in \Hilb^{(|\Index{I}| + 1) \times |\Index{J}|}$ and introduce an auxiliary matrix $\Matrix{T} = [\Matrix{V}_{\kappa}~\Matrix{e}_\omega - \truncate{\Batrix{G}}{\kappa}^\pinv \Batrix{r}] \in \F^{|\Index{J}| \times (\kappa + 1)}$, where $\omega$ is the local position of $j$ within $\Index{J}$.
This modification of $\Matrix{T}$ ensures that $\Batrix{K}$ has exactly the same block structure as in the second considered case.
An identical argument yields the same bound $\Norm{\Bentry{e}}{\Hilb} \leq \tfrac{\nu_{\kappa}(\Batrix{M})}{\nu_{\kappa}(\Batrix{G})} \sigma_{\kappa + 1}(\Batrix{A})$, where now
\begin{equation*}
    \det(\Matrix{T}^\ast \Matrix{T}) = \Norm{\Matrix{e}_\omega - \Matrix{V}_{\kappa} \Matrix{V}_{\kappa}^\ast \Matrix{e}_\omega}{\ell_2}^2 \leq 1.
\end{equation*}
Applying \cref{lemma:rank_excess} we obtain
\begin{equation*}
    \sum_{s = 1}^{|\Index{I}| + 1} \nu_{\kappa}^2(\Batrix{M}_{-s}) \geq \sum_{s = 1}^{|\Index{I}| + 1} \nu^2(\Batrix{M}_{-s} \Matrix{Q}_\kappa) = \Big(|\Index{I}| + 1 - \kappa + \Delta(\Batrix{M} \Matrix{Q}_\kappa) \Big) \nu_{\kappa}^2(\Batrix{M})
\end{equation*}
and can directly bound $\nu_{\kappa}(\Batrix{M}_{-s}) \leq \nu_{\kappa}(\Batrix{G})$ since $\Batrix{M}_{-s}$ is of size $|\Index{I}| \times |\Index{J}|$.
\end{proof}

A standard remark is in order that \Cref{theorem:cgr_error_maxvol} holds when the projective volume of $\Batrix{G}$ is maximum only among the submatrices that differ from it in at most one row and at most one column.

\subsection{Error bounds for cross approximation}
The preparatory work on the approximation guarantees for one-sided interpolation can now be streamlined to derive error bounds for cross approximation.
We begin with a bound in the $\ell_2(\Hilb)$ norm, the only unitarily invariant norm that is also invariant to transposition.

\begin{lemma}
\label{lemma:cross_l2_error}
Let $\Batrix{A} = \Batrix{B} + \Batrix{E} \in \Hilb^{m \times n}$ and suppose that
\begin{equation*}
    \Norm{\Batrix{B} - \Batrix{C} \Batrix{G}^\pinv \Batrix{R}}{\ell_2(\Hilb)} \leq \epsilon, \quad \Norm{\Batrix{B}^\trans - \Batrix{R}^\trans (\Batrix{G}^\trans)^\pinv \Batrix{C}^\trans}{\ell_2(\Hilb)} \leq \delta.
\end{equation*}
Then
\begin{equation*}
    \Norm{\Batrix{B} - \cross{\Batrix{A}}{\Index{I}, \Index{J}}}{\ell_2(\Hilb)} \leq \epsilon + \Norm{\Batrix{G}^\pinv \Batrix{R}}{\ell_2 \to \ell_2} (\delta + \Norm{\Batrix{C}_{\Batrix{E}}}{\ell_2(\Hilb)})
\end{equation*}
\end{lemma}
\begin{proof}
We expand the error by triangle inequality to obtain
\begin{align*}
    \Norm{\Batrix{B} - \cross{\Batrix{A}}{\Index{I}, \Index{J}}}{\ell_2(\Hilb)} &\leq \Norm{\Batrix{B} - \Batrix{C} \Batrix{G}^\pinv \Batrix{R}}{\ell_2(\Hilb)} + \Norm{\Batrix{C} \Batrix{G}^\pinv \Batrix{R} - \cross{\Batrix{A}}{\Index{I}, \Index{J}}}{\ell_2(\Hilb)} \\
    &\leq \epsilon + \Norm{\Batrix{C} - \big( (\Batrix{G}^\trans)^\pinv \Batrix{C}^\trans \big)^\trans \Batrix{G}}{\ell_2(\Hilb)} \Norm{\Batrix{G}^\pinv \Batrix{R}}{\ell_2 \to \ell_2} \\
    &= \epsilon + \Norm{\Batrix{C}^\trans - \Batrix{G}^\trans (\Batrix{G}^\trans)^\pinv \Batrix{C}^\trans}{\ell_2(\Hilb)} \Norm{\Batrix{G}^\pinv \Batrix{R}}{\ell_2 \to \ell_2} \\
    &\leq \epsilon + \big( \Norm{\Batrix{C}_{\Batrix{B}}^\trans - \Batrix{G}^\trans (\Batrix{G}^\trans)^\pinv \Batrix{C}^\trans}{\ell_2(\Hilb)} + \Norm{\Batrix{C}_{\Batrix{E}}}{\ell_2(\Hilb)} \big) \Norm{\Batrix{G}^\pinv \Batrix{R}}{\ell_2 \to \ell_2} \\
    &\leq \epsilon + \big( \Norm{\Batrix{B}^\trans - \Batrix{R}^\trans (\Batrix{G}^\trans)^\pinv \Batrix{C}^\trans}{\ell_2(\Hilb)} + \Norm{\Batrix{C}_{\Batrix{E}}}{\ell_2(\Hilb)} \big) \Norm{\Batrix{G}^\pinv \Batrix{R}}{\ell_2 \to \ell_2}.
\end{align*}
\end{proof}

The following theorem builds upon \cref{theorem:cgr_error}, and we simplify its error bound for the sake of presentation.
The proof relies on the sequential processing of a Bochner matrix and its transpose, so a symmetric bound holds too.
The effects induced by the regularization of two pseudoinverses can be incorporated in the cross-approximation bound by means of \cref{corollary:cgr_error_tau,corollary:cgr_error_k}, though we omit such extension for brevity.

\begin{theorem}
\label{theorem:cross_l2_error}
Let $\Batrix{A} = \Batrix{B} + \Batrix{E} \in \Hilb^{m \times n}$ with $\Batrix{B} \neq 0$, and let $\Batrix{B} = \Batrix{U} \Matrix{\Sigma} \Matrix{V}^\ast$ and $\Batrix{B}^\trans = \Batrix{P} \Matrix{\Lambda} \Matrix{Q}^\ast$ be SVDs.
Denote
\begin{gather*}
    \eta_r = \Norm{\Batrix{R}_{\Batrix{U}}^\pinv}{\ell_2(\Hilb) \to \ell_2}, \quad \eta_c = \Norm{\Matrix{C}_{\Matrix{V}^\ast}^\pinv}{\ell_2 \to \ell_2}, \quad \xi_r = \Norm{\Matrix{R}_{\Matrix{Q}}^\pinv}{\ell_2 \to \ell_2}, \quad \xi_c = \Norm{(\Batrix{C}_{\Batrix{P}^\trans}^\trans)^\pinv}{\ell_2(\Hilb) \to \ell_2}, \\
    \tilde{\eta} = \eta_r + \eta_c + 3 \eta_r \eta_c, \quad \tilde{\xi} = \xi_c + \xi_r + 3 \xi_c \xi_r.
\end{gather*}
Let $\Rank{T}{\Batrix{G}_{\Batrix{B}}} = \Rank{T}{\Batrix{B}}$, then
\begin{align*}
    \Norm{\Batrix{B} - \cross{\Batrix{A}}{\Index{I}, \Index{J}}}{\ell_2(\Hilb)} &\leq \tilde{\eta} (2 + \tilde{\xi}) \Norm{\Batrix{E}}{\ell_2(\Hilb)} \\
    &+ (1 + \tilde{\eta})(2 + \tilde{\xi}) \Big( \Norm{\Batrix{G}^\pinv}{\ell_2(\Hilb) \to \ell_2}  + \Norm{(\Batrix{G}^\trans)^\pinv}{\ell_2(\Hilb) \to \ell_2} \Big) \Norm{\Batrix{E}}{\ell_2(\Hilb)}^2.
\end{align*}
\end{theorem}
\begin{proof}
Coarsening of the bound in \cref{theorem:cgr_error} gives
\begin{equation*}
    \Norm{\Batrix{B} - \Batrix{C} \Batrix{G}^\pinv \Batrix{R}}{\ell_2(\Hilb)} \leq \tilde{\eta} \Norm{\Batrix{E}}{\ell_2(\Hilb)} + \Norm{\Batrix{G}^\pinv}{\ell_2(\Hilb) \to \ell_2} (1 + \tilde{\eta}) \Norm{\Batrix{E}}{\ell_2(\Hilb)}^2.
\end{equation*}
By the same token, we bound
\begin{equation*}
    \Norm{\Batrix{B}^\trans - \Batrix{R}^\trans (\Batrix{G}^\trans)^\pinv \Batrix{C}^\trans}{\ell_2(\Hilb)} \leq \tilde{\xi} \Norm{\Batrix{E}}{\ell_2(\Hilb)} + \Norm{(\Batrix{G}^\trans)^\pinv}{\ell_2(\Hilb) \to \ell_2} (1 + \tilde{\xi}) \Norm{\Batrix{E}}{\ell_2(\Hilb)}^2.
\end{equation*}
To finalize the proof via \cref{lemma:cross_l2_error}, we need a bound on $\Norm{\Batrix{G}^\pinv \Batrix{R}}{\ell_2 \to \ell_2}$, which we take from the end of the proof of \cref{theorem:cgr_error}:
\begin{equation*}
    \Norm{\Batrix{G}^\pinv \Batrix{R}}{\ell_2 \to \ell_2} \leq \eta_c + \Norm{\Batrix{G}^\pinv}{\ell_2(\Hilb) \to \ell_2} (1 + \eta_c) \Norm{\Batrix{E}}{\ell_2(\Hilb)}.
\end{equation*}
We substitute the bounds in \cref{lemma:cross_l2_error} and massage the expression using $\eta_c \leq \tilde{\eta}$.
\end{proof}

To extend the maximum-volume bound from \cref{theorem:cgr_error_maxvol} to cross approximation, we need to assume that the projective volumes of both $\Batrix{G}$ and $\Batrix{G}^\trans$ are simultaneously maximal across the submatrices of $\Batrix{A}$ and $\Batrix{A}^\trans$, respectively.
However, this is a stringent assumption, because submatrices maximizing both volumes at the same time do not exist for Bochner matrices in general.
Therefore, we rely instead on a \emph{quasi-maximum volume} assumption and require that $\theta \nu_\kappa(\Batrix{G})$ be larger than or equal to the maximum projective volume for $\theta \geq 1$.

\begin{theorem}
\label{theorem:cross_error_maxvol}
Let $\Batrix{A} \in \Hilb^{m \times n}$. 
Let $\rho,\kappa \in \N$ and $\theta_c, \theta_r \geq 1$.
Suppose $\Batrix{G}$ has the $\theta_c$-quasi-maximum $\kappa$-projective volume among all $|\Index{I}| \times |\Index{J}|$ submatrices of $\Batrix{A}$ and $\Batrix{G}^\trans$ has the $\theta_r$-quasi-maximum $\rho$-projective volume among all $|\Index{J}| \times |\Index{I}|$ submatrices of $\Batrix{C}^\trans$.
Denote $\Batrix{X} = (\truncate{\Batrix{G}^\trans}{\rho}^\pinv \Batrix{C}^\trans)^\trans \Batrix{G} (\truncate{\Batrix{G}}{\kappa}^\pinv \Batrix{R})$.
If $\nu_{\kappa}(\Batrix{G}) > 0$ and $\nu_{\rho}(\Batrix{G}^\trans) > 0$,
\begin{align*}
    \Norm{\Batrix{A} - \Batrix{X}}{\ell_{\infty}(\Hilb)} &\leq \theta_c \sqrt{1 + \frac{\kappa - \Delta_\kappa'}{|\Index{I}| + 1 - \kappa + \Delta_\kappa'}} \sqrt{1 + \frac{\kappa}{|\Index{J}| + 1 - \kappa}} \sigma_{\kappa + 1}(\Batrix{A}) \\
    &+ \theta_r \theta_c \sqrt{1 + \frac{\rho}{|\Index{I}| + 1 - \rho}} \sqrt{1 - \frac{1}{\theta_c^2} + \frac{\kappa}{|\Index{J}| + 1 - \kappa}} \sigma_{\rho + 1}(\Batrix{C}^\trans).
\end{align*}
with $\Delta_\kappa' = \min_{i \not\in \Index{I}} \min_{1 \leq j \leq n} \Delta_{\kappa}\big( \Batrix{A}(\Index{I} \cup \{i\}, \Index{J} \cup \{j\}) \big)$.
\end{theorem}
\begin{proof}
By triangle inequality, we have
\begin{equation*}
    \Norm{\Batrix{A} - \Batrix{X}}{\ell_{\infty}(\Hilb)} \leq \Norm{\Batrix{A} - \Batrix{C} \truncate{\Batrix{G}}{\kappa}^\pinv \Batrix{R}}{\ell_{\infty}(\Hilb)} + \Norm{\Batrix{C} \truncate{\Batrix{G}}{\kappa}^\pinv \Batrix{R} - \Batrix{X}}{\ell_{\infty}(\Hilb)}.
\end{equation*}
We bound the first term via a quasi-maximum-volume variant of \cref{theorem:cgr_error_maxvol}, the proof of which is a direct modification of the original proof:
\begin{equation*}
    \Norm{\Batrix{A} - \Batrix{C} \truncate{\Batrix{G}}{\kappa}^\pinv \Batrix{R}}{\ell_{\infty}(\Hilb)} \leq \theta_c \sqrt{1 + \frac{\kappa - \Delta_\kappa'}{|\Index{I}| + 1 - \kappa + \Delta_\kappa'}} \sqrt{1 + \frac{\kappa}{|\Index{J}| + 1 - \kappa}} \sigma_{\kappa + 1}(\Batrix{A}).
\end{equation*}
Suppose that the maximum error in the second term is attained at the $(i,j)$th entry and denote $\Batrix{c}^\trans = \Batrix{C}(i,:)$ and $\Batrix{r} = \Batrix{R}(:, j)$ as previously:
\begin{align*}
    \Norm{\Batrix{C} \truncate{\Batrix{G}}{\kappa}^\pinv \Batrix{R} - \Batrix{X}}{\ell_{\infty}(\Hilb)} &= \Norm{(\Batrix{C} - (\truncate{\Batrix{G}^\trans}{\rho}^\pinv \Batrix{C}^\trans)^\trans \Batrix{G}) \truncate{\Batrix{G}}{\kappa}^\pinv \Batrix{R}}{\ell_{\infty}(\Hilb)} \\
    &= \Norm{(\Batrix{c} - \Batrix{G}^\trans \truncate{\Batrix{G}^\trans}{\rho}^\pinv \Batrix{c})^\trans \truncate{\Batrix{G}}{\kappa}^\pinv \Batrix{r}}{\Hilb} \\
    &\leq \Norm{\Batrix{c} - \Batrix{G}^\trans \truncate{\Batrix{G}^\trans}{\rho}^\pinv \Batrix{c}}{\ell_2(\Hilb)} \Norm{\truncate{\Batrix{G}}{\kappa}^\pinv \Batrix{r}}{\ell_2}.
\end{align*}
As in the proof of \cref{theorem:cgr_error_maxvol}, we consider different cases based on the position of $(i,j)$ relative to $(\Index{I}, \Index{J})$.
When $i \in \Index{I}$ and $\Batrix{c} = \Batrix{G}^\trans \Matrix{e}_{\omega}$,
\begin{align*}
    \Norm{\Batrix{c} - \Batrix{G}^\trans \truncate{\Batrix{G}^\trans}{\rho}^\pinv \Batrix{c}}{\ell_2(\Hilb)} &= \Norm{(\Batrix{G}^\trans - \Batrix{G}^\trans \truncate{\Batrix{G}^\trans}{\rho}^\pinv \Batrix{G}^\trans) \Matrix{e}_{\omega}}{\ell_2(\Hilb)} \\
    &\leq \Norm{\Batrix{G}^\trans - \truncate{\Batrix{G}^\trans}{\rho}}{\ell_2 \to \ell_2(\Hilb)} = \sigma_{\rho + 1}(\Batrix{G}^\trans) \leq \sigma_{\rho + 1}(\Batrix{C}^\trans).
\end{align*}
When $i \not\in \Index{I}$, we repeat the argument from the third part of the proof of \cref{theorem:cgr_error_maxvol}.
Namely, we form an augmentation $\Batrix{M} = [\Batrix{G}^\trans~\Batrix{c}]$, take an auxiliary matrix $\Matrix{T}$, multiply $\Batrix{K} = \Batrix{M} \Matrix{T}$, obtain an equality
\begin{equation*}
    \nu(\Batrix{K}) = \nu_{\rho}(\Batrix{G}^{\trans}) \Norm{\Batrix{c} - \Batrix{G}^\trans \truncate{\Batrix{G}^\trans}{\rho}^\pinv \Batrix{c}}{\ell_2(\Hilb)},
\end{equation*}
and derive the bound
\begin{equation*}
   \Norm{\Batrix{c} - \Batrix{G}^\trans \truncate{\Batrix{G}^\trans}{\rho}^\pinv \Batrix{c}}{\ell_2(\Hilb)} = \frac{\nu(\Batrix{K})}{\nu_{\rho}(\Batrix{G}^{\trans})} \leq \frac{\nu_{\rho + 1}(\Batrix{M})}{\nu_{\rho}(\Batrix{G}^{\trans})} \leq  \theta_r \sqrt{1 + \frac{\rho}{|\Index{I}| + 1 - \rho}} \sigma_{\rho + 1}(\Batrix{C}^\trans).
\end{equation*}
When $j \in \Index{J}$ and $\Batrix{r} = \Batrix{G} \Matrix{e}_{\omega}$, we have
\begin{equation*}
    \Norm{\truncate{\Batrix{G}}{\kappa}^\pinv \Batrix{r}}{\ell_2} = \Norm{\truncate{\Batrix{G}}{\kappa}^\pinv \Batrix{G} \Matrix{e}_{\omega}}{\ell_2} \leq \Norm{\truncate{\Batrix{G}}{\kappa}^\pinv \Batrix{G}}{\ell_2 \to \ell_2} \leq 1.
\end{equation*}
When $j \not\in \Index{J}$, consider $\Batrix{M} = [\Batrix{G}~\Batrix{r}]$, and let $\truncate{\Batrix{G}}{\kappa} = \Batrix{U}_\kappa \Matrix{\Sigma}_\kappa \Matrix{V}_\kappa^\ast$ be an SVD.
Then, by the matrix determinant lemma,
\begin{align*}
    \nu_\kappa^2(\Batrix{M}) \geq \nu^2(\Batrix{U}_\kappa^\ast \Batrix{M}) &= \det\Big( \Matrix{\Sigma}_\kappa^2 + (\Batrix{U}_\kappa^\ast \Batrix{r}) (\Batrix{U}_\kappa^\ast \Batrix{r})^\ast \Big) \\
    &= \nu_\kappa^2(\Batrix{G}) (1 + \Norm{\Matrix{\Sigma}_\kappa^{-1} \Batrix{U}_\kappa^\ast \Batrix{r}}{\ell_2}^2) = \nu_\kappa^2(\Batrix{G}) (1 + \Norm{\truncate{\Batrix{G}}{\kappa}^\pinv \Batrix{r}}{\ell_2}^2),
\end{align*}
and consequently
\begin{equation*}
    \Norm{\truncate{\Batrix{G}}{\kappa}^\pinv \Batrix{r}}{\ell_2} \leq \sqrt{\frac{\nu_\kappa^2(\Batrix{M})}{\nu_\kappa^2(\Batrix{G})} - 1} \leq \sqrt{\theta_c^2 + \frac{\theta_c^2 \kappa}{|\Index{J}| + 1 - \kappa} - 1} = \theta_c \sqrt{1 - \frac{1}{\theta_c^2} + \frac{\kappa}{|\Index{J}| + 1 - \kappa}}.
\end{equation*}
\end{proof}

The bound in \cref{theorem:cross_error_maxvol} is asymmetric, and a similar one holds with the role of rows and columns reversed.
In addition, the statement still holds when $\Batrix{G}$ optimizes the volume among submatrices of $\Batrix{A}$ that differ from it in at most one row and at most one column, and $\Batrix{G}^\trans$ optimizes the volume among submatrices of $\Batrix{C}^\trans$ that differ from it in at most one column.
\section{Numerical experiments}
\label{sec:numerics}

\subsection{Adaptive Bochner cross approximation}
For classical matrices, \emph{adaptive cross approximation} (ACA, \cite{bebendorf2000approximation, bebendorf2003adaptive}) is a widely adopted approach to index selection for cross approximation, whereby columns and rows are picked one by one.
At its every iteration, the position of an entry of large absolute value is chosen in the current residual.
In full pivoting, the largest entry in the whole residual is selected, but then every entry needs to be inspected.
In practice, \emph{rook pivoting} delivers a good balance between stability and efficiency and is de facto the standard way to select indices in ACA \cite{dolgov2020parallel, shi2026distributed}.
\Cref{alg:rook} is a direct extension of rook pivoting to Bochner matrices.

\begin{algorithm2e}[ht]
\caption{Rook pivoting for Bochner matrices}
\label{alg:rook}
\DontPrintSemicolon
\KwIn{$\Batrix{A} \in \Hilb^{m \times n}$, starting column index $j_* \in [n]$, number of rounds $n_{\mathrm{rook}} \in \N_0$}
$i_* \gets \argmax_{i \in [m]} \Norm{\Batrix{A}(i, j_*)}{\Hilb}$\;
\For{$s = 1, \ldots, n_{\mathrm{rook}}$}{
    $j_* \gets \argmax_{j \in [n]} \Norm{\Batrix{A}(i_*, j)}{\Hilb}$\;
    $i_* \gets \argmax_{i \in [m]} \Norm{\Batrix{A}(i, j_*)}{\Hilb}$\;
}
\KwOut{pivot indices $(i_*, j_*)$}
\end{algorithm2e}

The ACA algorithm maintains a decomposition $\Matrix{A} = \Matrix{B} + \Matrix{E}$, starting with $\Matrix{B} = 0$ and increasing its rank by one at every iteration by adding a rank-one cross approximation of $\Matrix{E}$ to it.
A fundamental property of ACA is additivity; that is, if $\Matrix{B} = \cross{\Matrix{A}}{\Index{I}, \Index{J}}$ with invertible submatrix $\Matrix{A}(\Index{I}, \Index{J})$ and $\Matrix{E}(i,j) \neq 0$ then
\begin{equation*}
    \Matrix{B} + \cross{\Matrix{E}}{i,j} = \cross{\Matrix{A}}{\Index{I} \cup \{ i \}, \Index{J} \cup \{ j \}}.
\end{equation*}
Therefore, ACA iteratively builds up cross approximation based on rank-one updates.

However, this property ceases to hold for Bochner matrices.
Recall the example from \cref{sec:bmx_m} of a $2 \times 2$  Bochner matrix whose entries are mutually orthonormal:
\begin{equation*}
    \Batrix{B} = \cross{\Batrix{A}}{1,1} = \begin{bmatrix}
        \Batrix{A}(1,1) & 0 \\
        0 & 0
    \end{bmatrix}, \quad \Batrix{B} + \cross{\Batrix{E}}{2,2} = \begin{bmatrix}
        \Batrix{A}(1,1) & 0 \\
        0 & \Batrix{A}(2,2)
    \end{bmatrix},
\end{equation*}
whereas $\cross{\Batrix{A}}{\{1,2\}, \{1,2\}} = \Batrix{A}$.
Rank-one updates do not add up to a higher-rank cross approximation, requiring a proper recomputation.
We take this into account in our \emph{adaptive Bochner cross} (ABC) approximation algorithm (\Cref{alg:abc}).

\begin{algorithm2e}[ht]
\caption{Adaptive Bochner cross approximation}
\label{alg:abc}
\DontPrintSemicolon
\KwIn{$\Batrix{A} \in \Hilb^{m \times n}$, number of iterations $n_{\mathrm{iter}} \in \N$, number of sub-iterations $n_{\mathrm{sub}} \in \N$, number of rook-pivoting rounds $n_{\mathrm{rook}} \in \N_0$}
$\Index{I} \gets \varnothing,~\Index{J} \gets \varnothing,~\Batrix{B} \gets 0$\;
\For{$k = 1, \ldots, n_{\mathrm{iter}}$}{
    $\Batrix{E} \gets \Batrix{A} - \Batrix{B}$\tcc*{actual residual}
    \For{$s = 1, \ldots, n_{\mathrm{sub}}$}{
        $j_* \gets \mathrm{random}([n])$\;
        $(i_*,j_*) \gets \scup{rook}(\Batrix{E}, j_*, n_{\mathrm{rook}})$\;
        \If{$\Batrix{E}(i_*,j_*) \neq 0$}{
            $\Index{I} \gets \Index{I} \cup \{ i_* \},~\Index{J} \gets \Index{J} \cup \{ j_* \}$\;
            $\Batrix{\Delta} \gets \cross{\Batrix{E}}{i_*, j_*}$\tcc*{rank-one correction}
            $\Batrix{E} \gets \Batrix{E} - \Batrix{\Delta}$\tcc*{approximate residual}
        }
    }
    $\Batrix{B} \gets \cross{\Batrix{A}}{\Index{I}, \Index{J}}$\;
}
\KwOut{indices $(\Index{I}, \Index{J})$, cross approximation $\Batrix{B} = \cross{\Batrix{A}}{\Index{I}, \Index{J}}$}
\end{algorithm2e}

Note the presence of the inner loop, in which the residual undergoes fast approximate rank-one corrections.
When $n_{\mathrm{sub}} = 1$, the ABC algorithm extends an ACA-like algorithm developed in \cite{caiafa2010generalizing} from classical third-order tensors to Bochner matrices (neglecting one mode).
With $n_{\mathrm{sub}} > 1$, we attempt to garner more indices at reduced computational cost, making the algorithm more flexible.

Unlike classical ACA, the ABC algorithm can select the same indices multiple times.
This follows from the weaker interpolation properties of Bochner cross approximation (\Cref{theorem:cross_interpol}), rooted in the difference between column and row ranks.

\subsection{Reduced-order modeling: Nonlinear Stokes equations}
To validate the ABC algorithm as a non-intrusive ROM method for parametric PDEs, we consider the two-dimensional Stokes equations in $\Omega \subset \Real^2$, which describe the stationary flow of a fluid with dominating viscous forces.
Let $\Matrix{u} : \Omega \to \Real^2$ be the velocity and $p : \Omega \to \Real$ be the pressure.
When there are no external forces, the Stokes equations read
\begin{equation}
\label{eq:stokes}
    -\mathrm{div}(2 \nu \epsilon(\Matrix{u}) - p \Matrix{I}_2) = 0,~\mathrm{div}(\Matrix{u}) = 0~\text{ in }~\Omega,
\end{equation}
where $\epsilon(\Matrix{u}) = \tfrac{1}{2} \nabla \Matrix{u} + \tfrac{1}{2} (\nabla \Matrix{u})^\trans$ and $\nu$ is the viscosity, which we model as
\begin{equation*}
    \nu_{\alpha,\beta}(\Matrix{u}) = \nu_0 + \alpha \Norm{\nabla \Matrix{u}}{\ell_2}^{\beta},
\end{equation*}
making the equations nonlinear and parametric.
We partition the boundary as $\partial\Omega = \partial\Omega_D \cup \partial\Omega_N$ and impose the following boundary conditions:
\begin{equation}
\label{eq:stokes_bcs}
    \Matrix{u} = 0~\text{ on }~\partial\Omega_D, \qquad 
    -(2 \nu \epsilon(\Matrix{u}) - p \Matrix{I}_2) \Matrix{n} = p_0 \Matrix{n}~\text{ on }~\partial\Omega_N,
\end{equation}
where $\Matrix{n} : \partial\Omega_N \to \Real^2$ is the outward normal and $p_0 : \partial\Omega_N \to \Real$.

\begin{figure}[t]
\centering
\begin{subfigure}[b]{0.328\linewidth}
\centering
	\includegraphics[width=\linewidth]{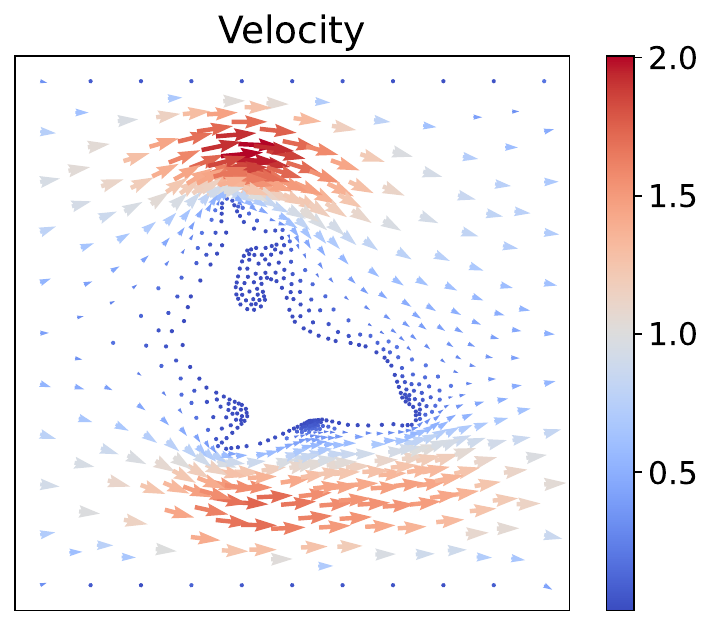}
\end{subfigure}%
\begin{subfigure}[b]{0.34\linewidth}
\centering
	\includegraphics[width=\linewidth]{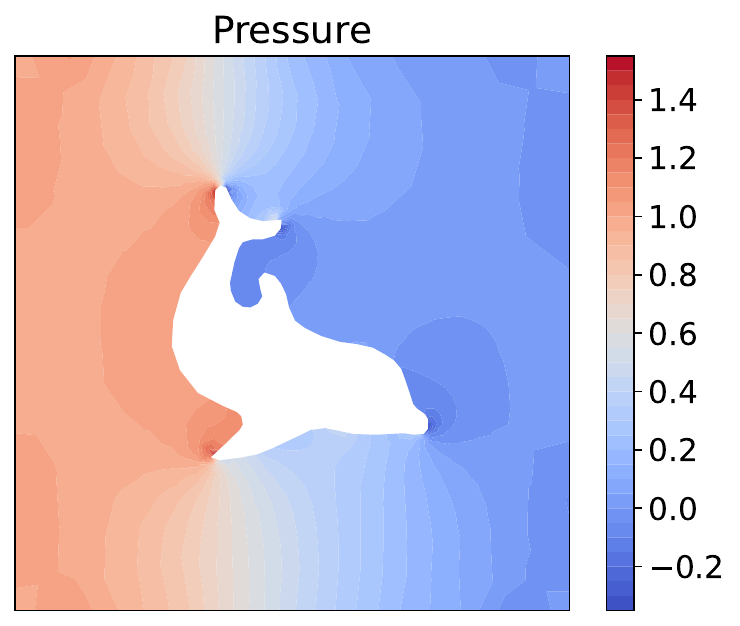}
\end{subfigure}%
\begin{subfigure}[b]{0.27\linewidth}
\centering
	\includegraphics[width=\linewidth]{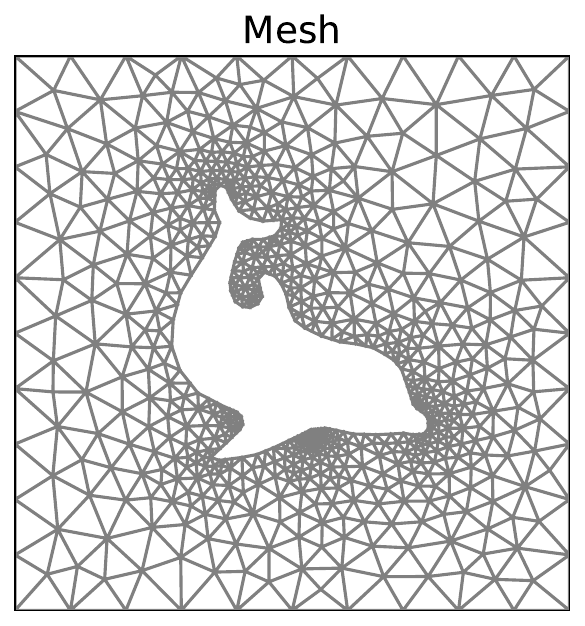}
\end{subfigure}
\caption{Finite-element solution to \cref{eq:stokes,eq:stokes_bcs} obtained with FEniCS as in \cite{rognes2017fenics}.}
\label{fig:fenics_solution}
\end{figure}

This boundary-value problem is a standard example used to illustrate the FEniCS\footnote{\url{https://fenicsproject.org/}} library.
We take the domain, mesh, and boundary conditions from \cite{rognes2017fenics} to make the problem \cref{eq:stokes,eq:stokes_bcs} concrete and use FEniCS to solve it using Taylor--Hood finite elements with 6495 degrees of freedom.
The mesh and the solution for $\nu_0 = 0.01$ and $\alpha = \beta = 0$ are shown in \cref{fig:fenics_solution}.

\begin{figure}[h]
\begin{subfigure}[b]{0.5\linewidth}
\centering
	\includegraphics[width=\linewidth]{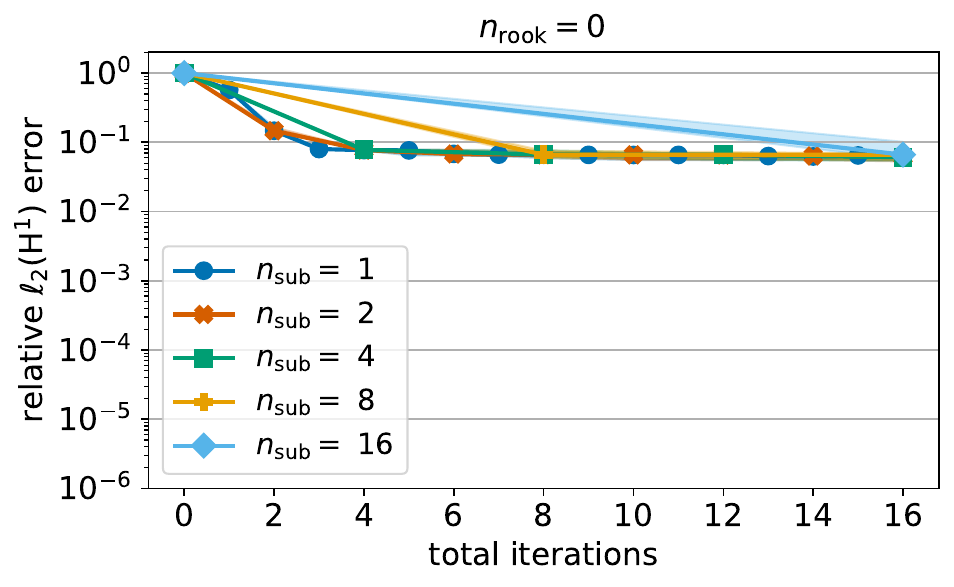}
\end{subfigure}\hfill
\begin{subfigure}[b]{0.48\linewidth}
\centering
	\includegraphics[width=\linewidth]{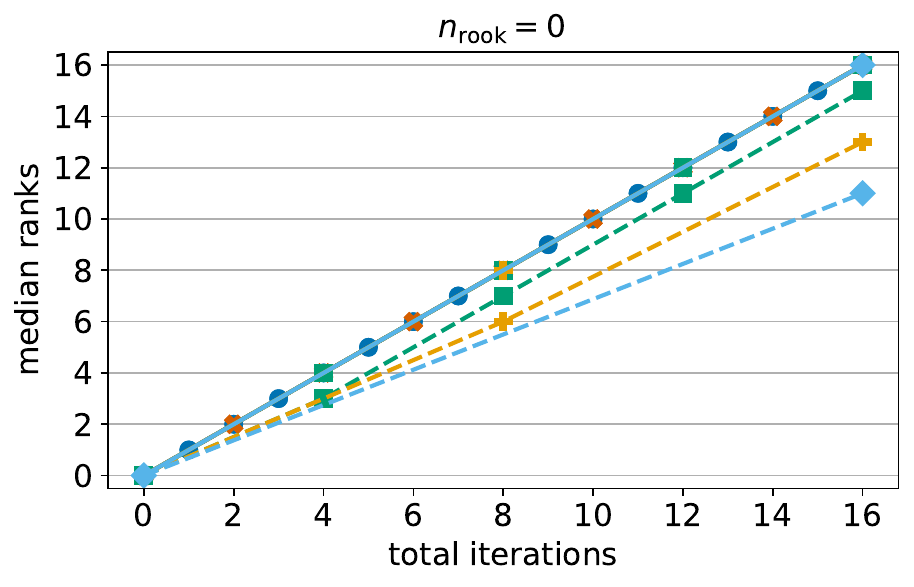}
\end{subfigure}
\begin{subfigure}[b]{0.5\linewidth}
\centering
	\includegraphics[width=\linewidth]{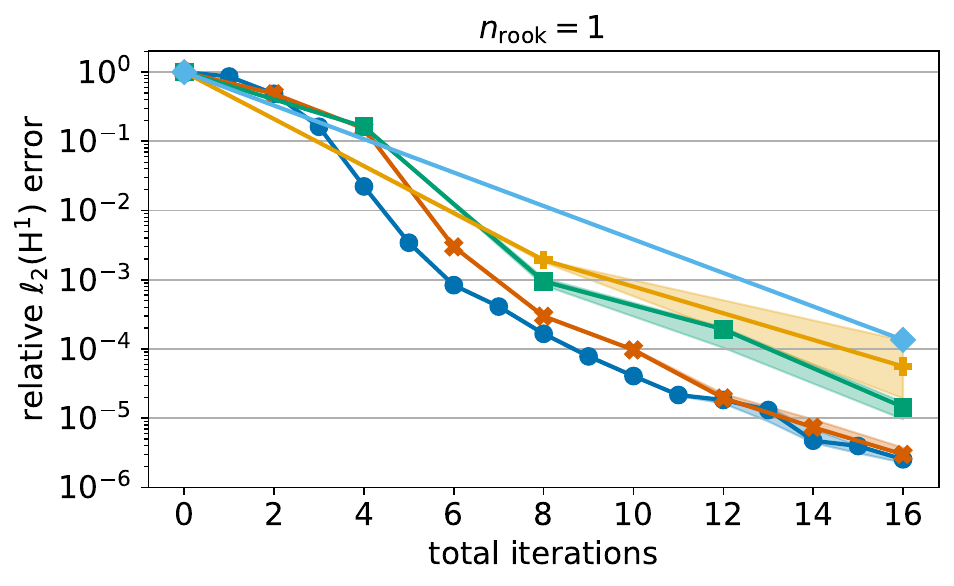}
\end{subfigure}\hfill
\begin{subfigure}[b]{0.48\linewidth}
\centering
	\includegraphics[width=\linewidth]{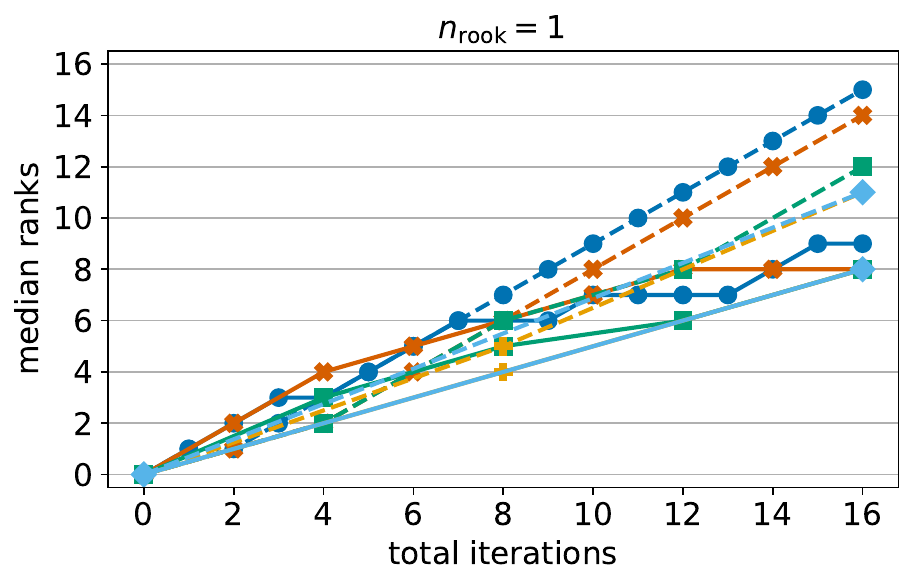}
\end{subfigure}
\caption{Cross approximation of the solution map of \cref{eq:stokes,eq:stokes_bcs} with ABC over 50 random seeds.
Left: median and 90th percentile of the relative approximation error.
Right: median of the row ranks (solid) and column ranks (dashed).
}
\label{fig:fenics_abc}
\end{figure}

We apply ABC to approximate the solution map of \cref{eq:stokes,eq:stokes_bcs} on a $100 \times 100$ uniform grid in the parameter space $(\alpha,\beta) \in [0,1] \times [0,1]$.
To define the corresponding Bochner matrix, we equip the 6495-dimensional space of finite-element solutions with the $\Hilb^1(\Omega; \Real^3)$ inner product.
In the experiments, we fix the product $n_{\mathrm{iter}} n_{\mathrm{sub}} = 16$ and vary $n_{\mathrm{sub}} \in \{ 1, 2, 4, 8, 16 \}$, plotting the statistics over 50 random seeds in \cref{fig:fenics_abc}.

When $n_{\mathrm{rook}} = 0$, the column selection is completely random and the approximation error stagnates, indicating that although the algorithm continues to blindly expand the index sets, the chosen subspaces fail to capture the structure of the solution map.
Conversely, when $n_{\mathrm{rook}} = 1$, both rows and columns are chosen adaptively, leading to decreasing approximation errors and fewer selected indices.
As expected, the errors decay faster per iteration for smaller values of $n_{\mathrm{sub}}$, because the algorithm relies less frequently on the approximate residual.
Finally, we observe that fewer rows are selected than columns, which signifies that the dependence on $\beta$ is stronger; indeed, $\beta$ controls the degree of nonlinearity.

To assess the accuracy of ABC approximation, we compare it with the quasioptimal HOSVD approximation of the same Tucker rank.
The results in \cref{tab:fenics_errors} for ABC with $n_{\mathrm{rook}} = 1$ and $n_{\mathrm{sub}} = 1$ show that it achieves errors that are about twice the HOSVD errors, despite querying the PDE solver only a small number of times.

\begin{table}[h]
    \centering
    \caption{Relative $\ell_2(\Hilb^1)$ approximation errors of the solution map of \cref{eq:stokes,eq:stokes_bcs} obtained with ABC and HOSVD for the same Tucker rank.}
    \label{tab:fenics_errors}
    \begin{tabular*}{\textwidth}{@{\extracolsep\fill}lcccc@{}}
        \toprule
        $n_{\mathrm{iter}}$ & 4 & 8 & 12 & 16 \\
        \midrule
        $(|\Index{I}|, |\Index{J}|)$ & $(3,3)$ & $(6,7)$ & $(7,11)$ & $(9,15)$ \\
        ABC & $2.2 \times 10^{-2}$ & $1.6 \times 10^{-4}$ & $1.8 \times 10^{-5}$ & $2.6 \times 10^{-6}$ \\
        HOSVD & $1.0 \times 10^{-2}$ & $7.3 \times 10^{-5}$ & $8.3 \times 10^{-6}$ & $1.3 \times 10^{-6}$ \\
        \bottomrule
    \end{tabular*}  
\end{table}
\section{Remarks}
\label{sec:remarks}

\emph{Unrestricted Hilbert space.}
The developed framework of Bochner matrices and their cross approximation does not impose any constraints on the underlying Hilbert space, relying exclusively on its inner-product structure.
In contrast, the framework of quasimatrices \cite{townsend2015continuous,shustin2022semi} and quasitensors \cite{larsen2024tensor,han2024guaranteed} requires a reproducing kernel Hilbert space or assumes that the entries are continuous functions.
Similarly, the infinite-dimensional extension of tubal tensors in \cite{mor2025quasitubal} requires separability.

\emph{Structured sampling.}
The formalism of Bochner matrices and tensors provides an abstraction for data models where the naturally accessible unit is a fiber or slice, rather than a single entry of a classical tensor.
Tensor completion with such structured sampling was studied in \cite{sorensen2025multilinear} and applied to tensorial parametric ROM of dynamical systems in \cite{mamonov2022interpolatory, mamonov2024slice, mamonov2024tensorial}.

\emph{Cross approximation of Bochner tensors.}
The present work establishes the theoretical foundation of cross approximation for Bochner matrices.
In a follow-up paper \cite{budzinskiy2025low}, we extend this framework to the Tucker-cross approximation of Bochner tensors, focusing primarily on its applications to ROM of parametric PDEs.

\section*{Acknowledgements}
I thank Vladimir Kazeev and Maxim Olshanskii for the discussions of tensorial parametric ROM, which have lead to this article.

\bibliographystyle{siamplain}
\bibliography{common/bib}

\end{document}